\documentclass[MSNbibl,number,citesort,dvips]{arxbj}
\usepackage{graphicx,url,breakurl}
\aid{0}
\volume{20}
\issue{3}
\pubyear{2014}
\firstpage{1484}
\lastpage{1506}
\doi{10.3150/13-BEJ530} 

\makeatletter
\newcommand{\RMo}{\mathrm{o}}

\newcommand{\mrmd}{\,\mathrm{d}}

\newtheorem{theorem}{Theorem}[section]

\newcommand{\QQ}{Q--Q}

\newcommand{\cal}{\mathcal}

\makeatother

\begin{document}
\begin{frontmatter}

\title{Comparison of multivariate distributions using
quantile--quantile plots and related tests}
\runtitle{Q--Q plots and related tests}

\begin{aug}
\author[1]{\inits{S.S.}\fnms{Subhra Sankar} \snm{Dhar}\corref{}\thanksref{1}\ead[label=e1]{dsubhra@gmail.com}},
\author[2]{\inits{B.}\fnms{Biman} \snm{Chakraborty}\thanksref{2}\ead[label=e2]{B.Chakraborty@bham.ac.uk}}\\ \and
\author[3]{\inits{P.}\fnms{Probal} \snm{Chaudhuri}\thanksref{3}\ead[label=e3]{probal@isical.ac.in}}
\runauthor{S.S. Dhar, B. Chakraborty and P. Chaudhuri} 
\address[1]{Department of Mathematics and Statistics,
Indian Institute of Technology, Kanpur 208016, India. \printead{e1}}
\address[2]{School of Mathematics, University of Birmingham, United
Kingdom.\\ \printead{e2}}
\address[3]{Theoretical Statistics and Mathematics Unit, Indian
Statistical Institute, Kolkata-700108, India. \printead{e3}}
\end{aug}

\received{\smonth{9} \syear{2011}}
\revised{\smonth{3} \syear{2013}}

%
\begin{abstract}
The univariate quantile--quantile (\mbox{\QQ}) plot is a well-known graphical
tool for examining whether two data sets are generated from the same
distribution or not. It is also used to determine how well a specified
probability distribution fits a given sample. In this article, we
develop and study a multivariate version of the \mbox{\QQ} plot based on the
spatial quantile. The usefulness of the proposed graphical device is
illustrated on different real and simulated data, some of which have
fairly large dimensions. We also develop certain statistical tests that
are related to the proposed multivariate \mbox{\QQ} plot and study their
asymptotic properties. The performance of those tests are compared with
that of some other well-known tests for multivariate distributions
available in the literature.
\end{abstract}

%

\end{frontmatter}

\section{Introduction}\label{sec1}
The univariate quantile--quantile (\mbox{\QQ}) plot is a diagnostic tool, which is
widely used to assess the distributional similarities and differences
between two independent samples (see, e.g., Gnanadesikan and Wilk \cite{WilGna68},
Gnanadesikan \cite{Gna77} and Chambers \textit{et al.} \cite{Chaetal83}).
As discussed in Doksum \cite{Dok74}, Doksum and Sievers \cite{DokSie76} and Koenker
(\cite{Koe05}, pages 31 and 32), there are some fundamental connections between the
\mbox{\QQ} plot and the two-sample problem involving a semi-parametric
treatment effect model. The \mbox{\QQ} plot is also a popular device for
checking the appropriateness of a specified probability distribution
for a given univariate data. While the univariate \mbox{\QQ} plot has a long
history as a graphical tool for data analysis, there are only limited
attempts in the literature to generalize the \mbox{\QQ} plot for multivariate
samples. One can construct the \mbox{\QQ} plot for multivariate
data using the marginal quantiles. However, a \mbox{\QQ} plot based on the marginal
quantiles fails to capture the nature of dependence among the marginals
of a
multivariate distribution. Such a \mbox{\QQ} plot can only compare the
marginal distributions, but it is inadequate for a proper comparison of two
multivariate distributions because the marginal quantiles do not characterize
a multivariate distribution (see the supplemental article (Dhar, Chakraborty and Chaudhuri~\cite{DhaChaCha}) for an illustrative example). 

Breckling and Chambers \cite{BreCha88}, Chaudhuri \cite{Cha96} and Koltchinskii \cite{Kol97}
extensively
studied a multivariate quantile, which is popularly known as the spatial
quantile. Koltchinskii (\cite{Kol97},
Corollary~2.9, page 446) established that these spatial quantiles
characterize multivariate distributions. In this article, we propose an
extension of the \mbox{\QQ} plot using the spatial quantiles
for multivariate data. As we will see in subsequent sections, these \mbox{\QQ}
plots are in many ways natural generalizations of the univariate \mbox{\QQ}
plot. In
particular, for a $d$-dimensional multivariate data, there will be
$d$ two-dimensional plots, where the points in each plot cluster around
a straight line with slope $= 1$ and intercept $= 0$ \textit{if and only
if} the two multivariate
distributions under comparison are identical.


Motivated by the one-sample \mbox{\QQ} plot, Shapiro and Wilk \cite{ShaWil65} proposed
a test for
normality of univariate data. We also propose and study some
statistical tests
for multivariate distributions, which are related to our multivariate
\mbox{\QQ} plots. In our numerical and asymptotic studies, those tests turn
out to have either comparable or superior performance when compared
with the Kolmogorov--Smirnov and the Cramer--von Mises tests for multivariate
distributions.

\section{Multivariate \mbox{\QQ} plots}\label{sec2}
Recall that a univariate \mbox{\QQ} plot based on two samples
with sizes $n$ and $m$ consists of $r$ ($r = n+m$ if $n \neq m$ and $r
= n$ if $n=m$) points in the two-dimensional plane, where for $i = 1,
2,\ldots, r$, the two coordinates of the $i$th point are the
$(i/r)$th quantiles of the two samples. Here, in order to compare the
quantiles, one has
to match the quantiles of one data set with the corresponding quantiles
of another
data set. Easton and McCulloch \cite{EasMcC90} made an attempt to solve a similar matching problem for
multivariate data.
Their procedure was based on the permutation of the data that produced
the minimum sum
of the Euclidean distances between the matching data points in the two
given samples.
Consequently, in order to assess how well a specified probability
distribution fits
a given multivariate sample, they used a sample simulated from the
specified distribution.
The \mbox{\QQ} plots proposed by them can be used in two-sample problems only if
the two samples have the same size. In this paper, we use a
matching procedure based on the spatial rank and the spatial quantile.
The procedure is computationally simple and can be used in a two-sample
problem even if the two samples do not have the same size. Further, in
the case of
a one-sample problem, where one tries to test whether a specified
distribution fits
the data well or not, the construction of our \mbox{\QQ} plot does not require
generation
of a sample from the specified distribution.

The spatial rank of ${\mathbf z}\in\mathbb{R}^{d}$ with respect to the
data cloud formed by the observations ${\cal{X}} = \{{\mathbf
x}_{1},\ldots, {\mathbf x}_{n}\}$ is defined as $n^{-1}\sum_{i:{\mathbf
x}_{i}\neq{\mathbf z}}\|{\mathbf z} - {\mathbf x}_{i}\|^{-1}({\mathbf
z} - {\mathbf x}_{i})$ (see, e.g., M{\"o}tt{\"o}nen and
Oja~\cite{MotOja95}, Chaudhuri \cite{Cha96} and Serfling \cite{Ser04}).
For a random vector ${\mathbf x}$ with a probability distribution $F$
on $\mathbb{R}^{d}$, the $d$-dimensional spatial quantile
$Q_{F}({\mathbf u}) = (Q_{F, 1}({\mathbf u}),\ldots, Q_{F, d}({\mathbf
u}))$ is defined as $Q_{F}({\mathbf u}) =  \arg\min_{Q\in\mathbb
{R}^{d}}E\{\Phi({\mathbf u}, {\mathbf x} - Q) - \Phi({\mathbf u},
{\mathbf x})\}$ (see Chaudhuri \cite{Cha96} and Koltchinskii
\cite{Kol97}). Here $\Phi({\mathbf u},{\mathbf s}) = \|{\mathbf s}\| +
\langle{\mathbf u}, {\mathbf s}\rangle$, ${\mathbf u}\in B^{d} =
\{{\mathbf v}\dvt {\mathbf v}\in\mathbb{R}^{d}, \|{\mathbf v}\|<1 \}$,
$\langle\cdot,\cdot\rangle$ is the Euclidean inner product, and
$\|\cdot\|$ is the Euclidean norm induced by the inner product. For a
random sample ${\cal{X}} = \{{\mathbf x}_1,\ldots, {\mathbf x}_n\}$,
the empirical spatial quantile $Q_{{\cal{X}}}({\mathbf u}) =
(Q_{{\cal{X}}, 1}({\mathbf u}),\ldots, Q_{{\cal{X}}, d}({\mathbf u}))$
is obtained by replacing $F$ with its empirical version $F_{n}$. When
different coordinate variables in a multivariate data are measured in
different units, the spatial quantiles and the spatial ranks are
usually computed after standardizing each coordinate variable
appropriately. Note\vadjust{\goodbreak} that when our objective is to compare the
distributions of two random vectors $\mathbf x$ and $\mathbf y$, the
problem is equivalent to comparing the distributions of ${\mathbf
A}^{-1}{\mathbf x}$ and ${\mathbf A}^{-1}{\mathbf y}$, where $\mathbf
A$ is any appropriate positive definite matrix used to standardize the
variables.

We now consider a one-sample multivariate problem involving a
$d$-dimensional data set
${\cal{X}} = {\{\mathbf x}_1,\ldots, {\mathbf x}_n\}$, where ${\mathbf x}_{i} =
(x_{i, 1},\ldots, x_{i, d})$ has distribution $F$, and let $F_{0}$ be
a specified probability distribution on $\mathbb{R}^{d}$. Let ${\mathbf
u}_1,\ldots, {\mathbf u}_n$
be the spatial ranks of the data points ${\mathbf x}_{i}, i = 1,\ldots,
n$. Suppose that $Q_{F_{0}}({\mathbf u}_k) = (Q_{F_{0}, 1}({\mathbf u}_k),\ldots, Q_{F_{0}, d}({\mathbf u}_k))$ is the ${\mathbf u}_k$th spatial
quantile of the specified distribution $F_{0}$,
where $k = 1,\ldots, n$. Note that since $Q_{{\cal{X}}}({\mathbf u}_{k}) =
{\mathbf x}_{k}$,
where $Q_{{\cal{X}}}({\mathbf u}_{k})$ is the ${\mathbf u}_{k}$th empirical
spatial quantile of the data set
${\cal{X}}$, a natural way of matching the quantiles of the data set
with those of the
specified probability distribution will be by setting the
correspondence between ${\mathbf x}_{k}$
and $Q_{F_{0}}({\mathbf u}_k)$ (see also Marden \cite{Mar98,Mar04}). Consider the
set of points in $\mathbb{R}^{2}$ defined as $S_{n, i}({\cal{X}},
F_{0}) = \{(x_{k, i}, Q_{F_{0}, i}({\mathbf u}_k))\dvt  k = 1,\ldots, n\}$,
where $Q_{F_{0}, i}({\mathbf u}_k)$ and $x_{k, i}$ are the $i$th
components of $Q_{F_{0}}({\mathbf u}_{k})$ and ${\mathbf x}_{k}$, respectively,
and $i = 1,\ldots, d$. In particular, when
$d = 1$, $S_{n, 1}({\cal{X}}, F_{0})$ coincides with the set of points
that form the univariate \mbox{\QQ} plot for the one-sample problem.
Theorem \ref{theor2.1}, stated
below, ensures that for all $i = 1,\ldots, d$, the points in the
$i$th two-dimensional plot will lie close to a straight line with
slope $= 1$ and intercept $= 0$ \textit{if
and only if} $F = F_{0}$.

\begin{theorem}\label{theor2.1}
Suppose that $F_{0}$ is a specified distribution having a positive
density function, which is bounded on every bounded subset of
$\mathbb{R}^{d}$ $(d\geq2)$, and the same is true for $F$, the true
distribution of the data. Assume that $S_{n, i}({\cal{X}}, F_{0})$ is
constructed using the ${\mathbf u}_{k}$'s lying in any given closed
ball in $\mathbb{R}^{d}$ with the center at the origin and the radius
strictly smaller than one. Let $L(\varepsilon)$ be the collection of
points that lie in an $\varepsilon$-neighborhood of a straight line with
slope $= 1$ and intercept $= 0$. Then, for every $\varepsilon> 0$, we have
\[
\lim_{n\rightarrow\infty}P \Biggl( \bigcap_{i = 1}^{d}
\bigl[S_{n,
i}({\cal{X}}, F_{0})\subseteq L(\varepsilon) \bigr]
\Biggr) = 1,
\]
if and only if $F = F_{0}$.
\end{theorem}

An implication of Theorem \ref{theor2.1} is that the plots constructed using
$S_{n, i}({\cal{X}}, F_{0})$ for $i = 1,\ldots, d$ can be used,
just like the univariate \mbox{\QQ} plot, to determine whether the
specified distribution $F_{0}$ fits the data well or not. In practice,
$F_{0}$ may involve some unspecified
parameters that need to be estimated from the data. For instance, there
may be some
unknown location and scatter parameters associated with $F_{0}$, and we can
estimate them using standard techniques like the maximum likelihood
method. In such a case, we can make an affine transformation of the data
using the maximum likelihood estimates of the location and the scatter
parameters. In view of the asymptotic consistency of the maximum likelihood
estimate under appropriate conditions, the assertion in Theorem \ref{theor2.1}
about the linearity of the \mbox{\QQ} plots remains
valid if we construct the \mbox{\QQ} plots using such transformed data, and
the data are actually
generated from $F_{0}$. One may also use other consistent
estimates of the location and the scale parameters having high
breakdown points (e.g.,
the minimum covariance determinant estimates; see Rousseeuw and Leroy \cite{RouLer87}), which are
robust against outliers. It will be appropriate to point out that
Easton and McCulloch \cite{EasMcC90}
also proposed an affine transformation of the data before constructing
their \mbox{\QQ} plots in the
one-sample problem. Their proposal is not related in any way to the
maximum likelihood
estimation based on the specified distribution $F_{0}$, and it involves
an iterative algorithm\vadjust{\goodbreak}
for computing the affine transformation. Easton and McCulloch \cite{EasMcC90}
did not consider the case
when the specified distribution involves unknown parameters other than
the location and the scatter parameters.
Any such parameter can be estimated by the maximum likelihood method
using $F_{0}$ and the data.

We next consider the two-sample multivariate problem involving two
independent $d$-dimensional data sets, namely, ${\cal{X}} = \{{\mathbf
x}_{1},\ldots, {\mathbf x}_{n}\}$ and ${\cal{Y}} = {\{\mathbf
y}_{1},\ldots, {\mathbf y}_{m}\}$, where ${\mathbf x}_{i} = (x_{i,
1},\ldots, x_{i, d})$ has distribution $F$, and ${\mathbf y}_{j} =
(y_{j, 1},\ldots, y_{j, d})$ has distribution $G$. Suppose that
${\mathbf u}_1,\ldots, {\mathbf u}_n$ and ${\mathbf u}_{n + 1},\ldots,
{\mathbf u}_{n + m}$ are the spatial ranks of these observations within
their respective data sets ${\cal{X}}$ and ${\cal{Y}}$, respectively.
As in the case of the one-sample problem, $Q_{{\cal{X}}}({\mathbf u}_k)
= {\mathbf x}_{k}$ for $k = 1,\ldots, n$, and $Q_{{\cal{Y}}}({\mathbf
u}_k) = {\mathbf y}_{k}$ for $k = n + 1,\ldots, n + m$. We compute
$Q_{{\cal{X}}}({\mathbf u}_k)$ for $k = n + 1,\ldots, n + m$ and
$Q_{{\cal{Y}}}({\mathbf u}_{k})$ for $k = 1,\ldots, n$ using the
algorithm given in Chaudhuri (\cite{Cha96}, pages 864 and~865). Then,
we can match the two sets of quantiles by setting the correspondence
between $Q_{{\cal{X}}}({\mathbf u}_k)$ and $Q_{{\cal{Y}}}({\mathbf
u}_k)$ for $k = 1,\ldots, n + m$. As in the case of the one-sample
problem, one may construct the \mbox{\QQ} plots for the two-sample
problem as a collection of $d$ two-dimensional plots, where each plot
corresponds to a component of the spatial quantile. Let $S_{n, m,
i}({\cal{X}}, {\cal{Y}}) = \{ (Q_{{\cal{X}}, i}({\mathbf u}_k),
Q_{{\cal{Y}}, i}({\mathbf u}_k))\dvt  k = 1,\ldots, (n + m)\}$, where
$Q_{{\cal{X}}, i}({\mathbf u}_k)$ and $Q_{{\cal {Y}}, i}({\mathbf
u}_k)$ are the $i$th components of $Q_{{\cal{X}}}({\mathbf u}_k)$ and
$Q_{{\cal{Y}}}({\mathbf u}_k)$, respectively, and $i = 1,\ldots, d$.
Note that when $d = 1$, our proposed multivariate matching coincides
with the usual way of matching the univariate quantiles in a two-sample
problem, and the points in $S_{n, m, 1}({\cal{X}}, {\cal{Y}})$ are same
as those used in constructing the univariate two-sample \mbox{\QQ}
plot. Theorem \ref{theor2.2}, stated below, ensures that for all $i =
1,\ldots, d$, the points in the $i$th two-dimensional plot will lie
close to a straight line with slope $= 1$ and intercept $= 0$
\textit{if and only if} $F = G$.

\begin{theorem}\label{theor2.2}
Suppose that $F$ and $G$ have positive density
functions, which are bounded on every bounded subset of $\mathbb
{R}^{d}$ $(d\geq2)$, and $S_{n, m, i}({\cal{X}}, {\cal{Y}})$ is
constructed using the ${\mathbf u}_{k}$'s lying in
any given closed ball in $\mathbb{R}^{d}$ with the center at the origin
and the radius strictly smaller than one. Further, let
$L(\varepsilon)$ be the collection of points that lie in an
$\varepsilon$-neighborhood of a straight line with slope $= 1$ and
intercept $= 0$, and assume that
$n, m\rightarrow\infty$ in such a way that $ \lim_{n,
m\rightarrow\infty}\frac{n}{(n + m)} = \lambda\in(0, 1)$.
Then, for every $\varepsilon> 0$, we have
\[
\lim_{n, m\rightarrow\infty}P \Biggl( \bigcap_{i = 1}^{d}
\bigl[S_{n, m,
i}({\cal{X}}, {\cal{Y}})\subseteq L(\varepsilon) \bigr] \Biggr)
= 1,
\]
if and only if $F = G$.
\end{theorem}

In view of the equivariance of the spatial quantiles under location and
homogeneous scale transformations, the assertions in Theorems \ref{theor2.1} and
\ref{theor2.2} will also hold for the straight line with slope $= \sigma$ and
intercept $= \mu_{i}$ ($i = 1,\ldots, d$) \emph{if and only if} $F({\mathbf
x}) = F_{0}(({\mathbf x} - \bolds\mu)/\sigma)$
and $F({\mathbf x}) = G(({\mathbf x} - \bolds\mu)/\sigma)$, respectively,
where $\bolds\mu = (\mu_1,\ldots, \mu_d)\in\mathbb
{R}^{d}$ and $\sigma> 0$.


We now briefly discuss some earlier attempts to develop graphical tools
for comparing multivariate
distributions. For bivariate data, Marden \cite{Mar98,Mar04} proposed a
version of the \mbox{\QQ} plot, which is
based on drawing arrows from the spatial quantiles in one sample to the
corresponding spatial
quantiles in another sample in a two-sample problem (or to the
corresponding spatial quantiles of
a specified probability distribution in a one-sample problem). However,
such an arrow plot can be
drawn only for a bivariate data. Also, when the two samples are related
to each other by a location and a
homogeneous scale transformation, such arrow plots cannot detect that
unlike our \mbox{\QQ}
plots. Friedman and Rafsky \cite{FriRaf81} proposed a different visualization
procedure for comparing the distributions of two multivariate samples.
Their methodology is based on the idea of a minimal spanning
tree. Liu, Parelius and Singh \cite{LiuParSin99} proposed an alternative
visualization device called the DD-plot
for comparing two multivariate data sets based on the concept of data
depth. However, none
of these graphical tools developed by Marden \cite{Mar98,Mar04}, Friedman and Rafsky \cite{FriRaf81} and
Liu, Parelius and Singh \cite{LiuParSin99} will coincide with the usual univariate \mbox{\QQ} plot when they
are applied to the univariate data,
and none of them can be taken as a natural multivariate extension of
the univariate \mbox{\QQ} plot.

\section{Tests for comparing multivariate distributions}\label{sec3}
For each two-dimensional plot in our \mbox{\QQ} plots, the overall deviation
of the points from the straight
line with slope $= 1$ and intercept $= 0$ can be measured by $\int\{
Q_{{\cal{X}}, i}({\mathbf u}) - Q_{F_{0}, i}({\mathbf u})\}^{2}\mrmd {\mathbf u}$
and $\int\{Q_{{\cal{X}}, i}({\mathbf u}) - Q_{{\cal{Y}}, i}({\mathbf u})\}
^{2}\mrmd {\mathbf u}$ for the one-sample and the two-sample
problems, respectively, where $i = 1,\ldots, d$. These deviations in
$d$ different plots
can be aggregated as $\sum_{i = 1}^{d}\int\{Q_{{\cal{X}},
i}({\mathbf u}) - Q_{F_{0}, i}({\mathbf u})\}^{2}\mrmd {\mathbf u} =
\int\|Q_{\cal{X}}({\mathbf u}) - Q_{F_{0}} ({\mathbf u})\|^{2}\mrmd {\mathbf u}$ and
$\sum_{i = 1}^{d}\int\{Q_{{\cal{X}}, i}({\mathbf u}) - Q_{{\cal{Y}},
i}({\mathbf u})\}^{2}\mrmd {\mathbf u} =
\int\|Q_{{\cal{X}}}({\mathbf u}) - Q_{{\cal{Y}}}({\mathbf u})\|^{2}\mrmd {\mathbf u}$
for the one-sample and the
two-sample problems, respectively. These aggregated quantities can be
taken as the total
deviations in our \mbox{\QQ} plots. These measures of total deviations can be
used to construct tests for
comparing multivariate distributions. Such tests will be rotationally
invariant in view of the rotational equivariance of
the spatial quantiles.

Let ${\cal{X}} = \{{\mathbf x}_1,\ldots, {\mathbf x}_n\}$ consist of i.i.d.
observations from an unknown
distribution $F$ having a density function, which is assumed to be
bounded on every bounded subset of $\mathbb{R}^{d}$ ($d \geq2$).
Suppose that we want to test $H_0\dvt  F = F_{0}(\Leftrightarrow
Q_{F}({\mathbf u}) = Q_{F_{0}} ({\mathbf u})$
for all ${\mathbf u}\in B^{d}$) against the\vadjust{\goodbreak}
alternative $H_1\dvt  F\neq F_{0}(\Leftrightarrow Q_{F}({\mathbf u})\neq
Q_{F_{0}} ({\mathbf u})$ for some ${\mathbf u}\in B^{d}$), where $F_{0}$
is a specified distribution having a density function, which is bounded
on every bounded subset of $\mathbb{R}^{d}$ ($d \geq2$). In
order to test $H_{0}$ against $H_{1}$, we can use the test statistic
$V_{n} = n\int\|Q_{\cal{X}}({\mathbf u}) - Q_{F_{0}} ({\mathbf u})\|^{2}\mrmd {\mathbf
u}$, where the
integral is over a closed ball with the center at the origin and the
radius strictly smaller than one.
Note that the test statistic $V_{n}$ (as well as the test statistic
$T_{n, m}$ considered later in this section) can be viewed as the sum
of the
arrow lengths in the arrow plot considered by Marden \cite{Mar98} for a
bivariate data.

Consider now a multivariate
Gaussian process $Z_{1}({\mathbf u})$ having zero mean and the covariance kernel
\[
k_{1}({\mathbf u}_1, {\mathbf u}_2) =
\bigl[D_{1}\bigl\{Q_{F_{0}}({\mathbf u}_1)\bigr\}
\bigr]^{-1}\bigl[D_{2}\bigl\{Q_{F_{0}}({\mathbf
u}_1), Q_{F_{0}}({\mathbf u}_2), {\mathbf
u}_1, {\mathbf u}_2\bigr\}\bigr] \bigl[D_{1}
\bigl\{Q_{F_{0}}({\mathbf u}_2)\bigr\}\bigr]^{-1}.
\]
%
Here $D_{1}\{Q_{F_{0}}({\mathbf u})\} = E_{F_{0}}[\|{\mathbf x}- Q_{F_{0}}({\mathbf
u})\|^{-1}\{I_{d} - \|{\mathbf x} - Q_{F_{0}}({\mathbf u})\|^{-2}({\mathbf x} -
Q_{F_{0}}({\mathbf u}))({\mathbf x} -\break  Q_{F_{0}}({\mathbf u}))^{T}\}]$, $D_{2}\{
Q_{F_{0}}({\mathbf u}), Q_{F_{0}}({\mathbf v}), {\mathbf u}, {\mathbf v}\} =
E_{F_{0}}[\{\|{\mathbf x} - Q_{F_{0}}({\mathbf u})\|^{-1}({\mathbf x} -
Q_{F_{0}}({\mathbf u})) + {\mathbf u}\}\{\|{\mathbf x} - Q_{F_{0}}({\mathbf
v})\|^{-1}({\mathbf x} - Q_{F_{0}}({\mathbf v})) + {\mathbf v}\}^{T}]$. Henceforth,
$I_{d}$ denotes the $d\times d$ identity matrix, all vectors are
assumed to be column vectors, and the superscript $T$ denotes the
transpose of a vector.
Let ${\cal{V}} = \int\|Z_{1}({\mathbf u})\|^{2}\mrmd {\mathbf u}$, where the
integral is over the same closed ball as in the definition of $V_{n}$.
We now state a theorem describing the asymptotic behaviour of the test
based on $V_{n}$.


\begin{theorem}\label{theor3.1}
Let $c_{1}(\alpha)$ be the $(1 - \alpha)$th quantile $(0 < \alpha<
1)$ of the distribution of ${\cal{V}}$.
A test, which rejects $H_{0}$ for $V_{n} > c_{1}(\alpha)$, will have
asymptotic size $\alpha$. Further,
when $H_{1}$ is true, the asymptotic power of the test will be one if
the integral defining $V_n$ is taken
over an appropriately large closed ball in $\mathbb{R}^{d}$.
\end{theorem}

In order to implement our test, we need to compute $V_{n}$, and we have
approximated the integral
that appears in this test statistic by an average of the integrand over
$1000$ i.i.d. Monte Carlo replications
obtained from the random generations of ${\mathbf u}$ from the uniform
distribution on a closed ball with the center
at the origin and the radius $= 0.99$. In view of the asymptotic
Gaussian distribution of the process
$\sqrt{n}\{Q_{{\cal{X}}}({\mathbf u}) - Q_{F_{0}}({\mathbf u})\}$ under $H_{0}$
and the well-known orthogonal
decomposition of a finite-dimensional
multivariate normal distribution, the distribution of the test
statistic $V_{n}$ under $H_{0}\dvt  F = F_{0}$
can be approximated by a weighted sum of chi-square random variables
each with one degree of freedom.
In our numerical work, we have computed $c_{1}(\alpha)$ by generating
1000 Monte Carlo replications
from a weighted sum of chi-square variables, where the weights are the
eigenvalues of the covariance
matrices of appropriate normal random vectors. Note that the covariance matrices
involve the spatial quantiles and certain expectations under the
specified distribution $F_{0}$, and those
can be computed numerically. 


Let us next consider a two-sample problem with two independent sets of
i.i.d. observations
${\cal{X}} = \{{\mathbf x}_1,\ldots, {\mathbf x}_n\}$ and
${\cal{Y}} = \{{\mathbf y}_1,\ldots, {\mathbf y}_m\}$ from the distributions
$F$ and $G$, respectively. We assume the same conditions on the density
functions of $F$ and $G$ as for the density functions of $F$ and
$F_{0}$ in the one-sample problem discussed above. In this two-sample
problem, our hypotheses are
$H_{0}^{*}\dvt  F = G (\Leftrightarrow Q_{F}({\mathbf u}) = Q_{G}({\mathbf u})$
for all ${\mathbf u}\in B^{d}$) and $H_{1}^{*}\dvt  F\neq G (\Leftrightarrow
Q_{F}({\mathbf u})\neq Q_{G}({\mathbf u})$
for some ${\mathbf u}\in B^{d}$). In order to test $H_{0}^{*}$ against
$H_{1}^{*}$, one
can use the test statistic $T_{n, m} = (n + m)\int\|Q_{{\cal{X}}}({\mathbf
u}) - Q_{{\cal{Y}}}({\mathbf u})\|^{2}\mrmd {\mathbf u}$, where the integral is over
a closed ball with the center at the origin and the radius strictly
smaller than one.

Let $Z_{2}({\mathbf u})$ be a multivariate Gaussian process having zero
mean and the covariance kernel
\[
k_{2}({\mathbf u}_1, {\mathbf u}_2) =
\frac{[D_{1}\{Q_{F}({\mathbf u}_1)\}
]^{-1}[D_{2}\{Q_{F}({\mathbf u}_1), Q_{F}({\mathbf u}_2), {\mathbf u}_1, {\mathbf u}_2\}
][D_{1}\{Q_{F}({\mathbf u}_2)\}]^{-1}}{\lambda(1 - \lambda)},
\]
%
where $\lambda$ is as defined in the statement of Theorem \ref{theor2.2}, and
$D_1$, $D_2$ are as defined before the statement of Theorem \ref{theor3.1}. Define
${\cal{T}} = \int\|Z_{2}({\mathbf u})\|^{2}\mrmd {\mathbf u}$, where the integral
is over the same closed ball as in the definition of $T_{n, m}$. We now
state a theorem describing the asymptotic behaviour of the test based
on $T_{n, m}$.


\begin{theorem}\label{theor3.2}
Let $c_{2}(\alpha)$ be the $(1 - \alpha)$th quantile $(0 < \alpha<
1)$ of the distribution of $\cal{T}$.
A test, which rejects $H_{0}^{*}$ for $T_{n, m} > c_{2}(\alpha)$, will
have asymptotic size $\alpha$. Further, when $H_{1}^{*}$ is true, the
test will have asymptotic power one if the integral defining $T_{n, m}$
is taken over an appropriately large closed ball in $\mathbb{R}^{d}$.
\end{theorem}

For numerical implementation, one can compute $T_{n, m}$ and
$c_{2}(\alpha)$
for the two-sample problem in a similar way as we have computed $V_{n}$
and $c_{1}(\alpha)$, respectively,
in the one-sample problem. However, here we have estimated the unknown
quantities (i.e., the spatial quantiles and certain expectations under
$H_{0}^{*}$) appearing in the
covariance kernel based on the combined sample of the ${\mathbf x}$'s and
the ${\mathbf y}$'s.
In\vadjust{\goodbreak} Sections~\ref{sec5} and \ref{sec6}, we have compared the performance of our tests
with that of the Kolmogorov--Smirnov and the Cramer--von Mises tests for
multivariate
distributions. For numerical implementation, we have used $R$ codes
that are available from the
first author of the paper.

\section{Demonstration of multivariate \mbox{\QQ} plots using simulated and
real data}\label{sec4} 
\begin{figure}

\includegraphics{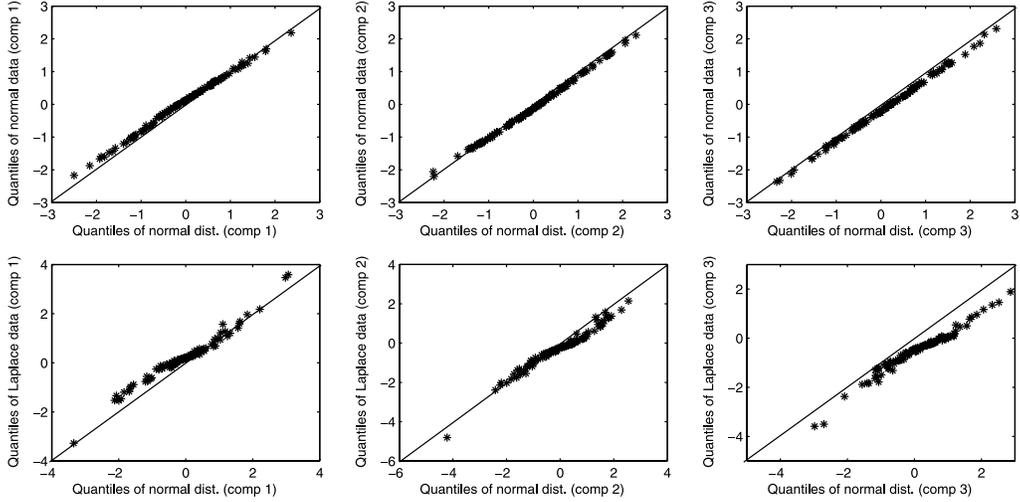}

\caption{The \mbox{\QQ} plots for the one-sample examples, where the specified
distribution is trivariate normal. The plots in the first and the
second rows are for the examples, where the distributions of the data
are trivariate
normal and trivariate Laplace, respectively.}\label{fig1}
\end{figure}

We begin with the one-sample problem and
consider two simulated data sets each consisting of 100 i.i.d.
observations. The observations in the first set were generated from the
trivariate normal distribution having zero mean and scatter matrix
$\Sigma= ((\sigma_{ij}))_{1\leq i, j\leq3}$ with $\sigma_{i, i} = 1$
for $i = 1, 2, 3$, $\sigma_{1, 2} = 0.5$, $\sigma_{1, 3} = 0.2$ and
$\sigma_{2, 3} = 0.3$. For the second set, the observations were
generated from the trivariate Laplace distribution with p.d.f.
$f({\mathbf x}) = (1/8\pi)\exp^{-\|{\mathbf x}\|}$. For both of them,
we considered the trivariate normal distribution as the specified
distribution $F_0$ with unknown parameters $\bolds\mu$ and $\Sigma$.
Following the remarks after Theorem \ref{theor2.1}, $\bolds\mu$ and
$\Sigma$ were estimated from each data set using the sample mean vector
and the sample dispersion matrix, respectively, which are the maximum
likelihood estimates in this case. We standardized the data sets using
these estimates and compared the spatial quantiles of the standardized
data with those of the standard trivariate normal distribution. We
computed the spatial quantiles for standard trivariate normal
distributions using the results in Marden (\cite{Mar98}, pages 824 and
825). The \mbox{\QQ} plots for the two simulated data sets are
displayed in Figure~\ref{fig1}. 

It is clearly evident from the plots in the first row of Figure~\ref{fig1} that the
specified distribution fits the data well
as the points in those plots are tightly clustered around the
straight\vadjust{\goodbreak}
line with slope $= 1$ and intercept $= 0$. On the other hand, in each
\mbox{\QQ} plot in the
second row, the points are significantly deviating from the straight
line with slope $= 1$ and intercept $= 0$, and the points are actually clustered
around a nonlinear curve. 
We have also computed the $p$-values for the one-sample test discussed in
Section~\ref{sec3} for testing $H_{0}\dvt  F = F_0$ against $H_1\dvt  F\neq F_0$ for
these two simulated data sets. We have obtained a high $p$-value $=
0.784$ for the first sample whereas the $p$-value for the second
example is $0.049$, which is quite small. 

\begin{figure}

\includegraphics{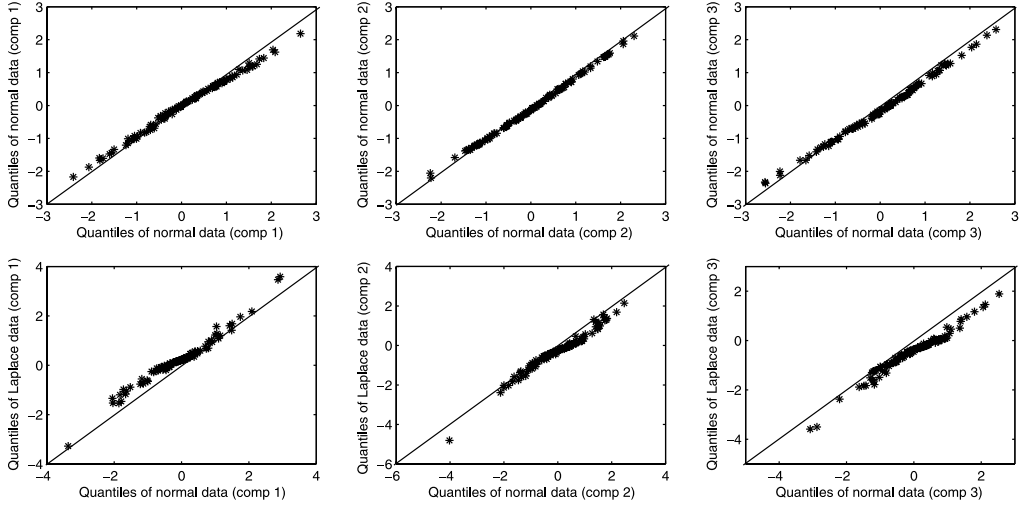}

\caption{The \mbox{\QQ} plots for the two-sample problem. The plots in the
first row for an example, where the samples are generated from the same
distribution, and those in the second row for an example, where the
samples are generated from different distributions.}\label{fig2}
\end{figure}

We next consider two simulated data sets to demonstrate our \mbox{\QQ} plots
for the two-sample problem. In both the data sets, the distribution of
the first sample $F$ was chosen to be the standard trivariate normal
distribution while $G$, the distribution of the second sample, was
taken to be the standard trivariate normal in one set and the
trivariate Laplace distribution in the other set. The size of each
sample was 100. 
The \mbox{\QQ} plots for the two data sets are displayed in Figure~\ref{fig2}. In each
plot in the
first row of Figure~\ref{fig2}, the points are tightly clustered around the
straight line with slope $= 1$ and intercept $= 0$. On the other hand,
the points are significantly deviating from the straight line with
slope $= 1$ and intercept $= 0$ in each plot in the second row of
Figure~\ref{fig2}.
We also carried out the two-sample test described in Section~\ref{sec3} for
testing $H_{0}^{*}\dvt  F = G$ against $H_{1}^{*}\dvt  F\neq G$, and we
obtained a high $p$-value $= 0.731$ for the first data set whereas a
small $p$-value $= 0.048$ was obtained for the second data set.

\begin{figure}[b]

\includegraphics{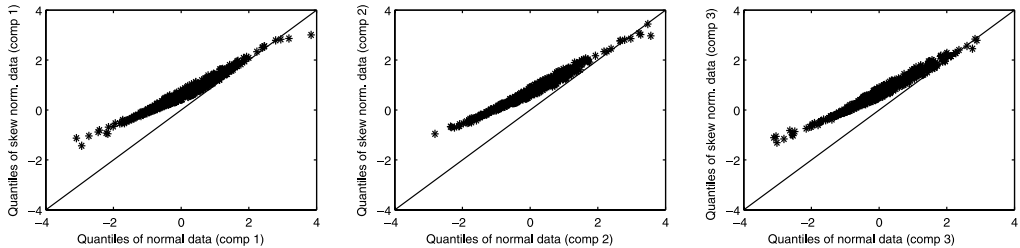}

\caption{The \mbox{\QQ} plots for the two-sample problem, where the first
sample is generated from the standard trivariate
normal distribution, and the second sample is generated from a
trivariate skew-normal
distribution.}\label{fig3}
\end{figure}

\subsection{Detection of special features using multivariate \mbox{\QQ} plots}\label{sec4.1}
We now consider a two-sample problem, where the
first sample consists of 100 i.i.d. observations from the standard
trivariate normal distribution ($F$),
and the second sample consists of 100 i.i.d. observations from a
trivariate skew-normal distribution ($G$)
(see Azzalini and Dalla Valle (\cite{AzzDal96}, page 717)). The p.d.f. of the
trivariate skew-normal
distribution is given by $f({\mathbf z}) = 2\phi_{3}({\mathbf z}; \Omega)\Phi
(\alpha^{T}{\mathbf z})$, where ${\mathbf z}\in\mathbb{R}^{3}$,
$\alpha^{T} = \frac{\lambda^{T}\Psi^{-1}\Delta^{-1}}{\sqrt{1 + \lambda
^{T}\Psi^{-1}\lambda}}$,
$\Delta= \operatorname{diag} (\sqrt{1 - \delta_{1}^{2}}, \sqrt{1 - \delta
_{2}^{2}}, \sqrt{1 - \delta_{3}^{2}} )$,
$\lambda=  (\frac{\delta_1}{\sqrt{1 - \delta_1^{2}}}, \frac{\delta
_2}{\sqrt{1 - \delta_2^{2}}}, \frac{\delta_3}{\sqrt{1 - \delta
_3^{2}}} )^{T}$, and
$\Omega= \Delta(\Psi+ \lambda\lambda^{T})\Delta$. Here $\phi_{3}({\mathbf
z}; \Omega)$ denotes the
p.d.f. of a trivariate normal distribution with standardized marginals
and correlation matrix
$\Omega$, and $\Phi$ is the distribution function of the standard
univariate normal distribution.
In this study, we have considered $\delta_{1} = \delta_{2} = \delta_{3}
= 0.9$ and $\Psi= I_{d}$. The \mbox{\QQ} plots for this two-sample problem
are displayed in Figure~\ref{fig3}, and we see a heavier tail in one direction
in each plot in this figure. This is an indication that one sample is
generated from a more skewed distribution than the other. Also,
the small $p$-value $= 0.048$ obtained using our two-sample test for
testing $H_{0}^{*}\dvt  F = G$ against $H_{1}^{*}\dvt  F \neq G$ implies that the
two distributions are significantly different in this data set.

\begin{figure}

\includegraphics{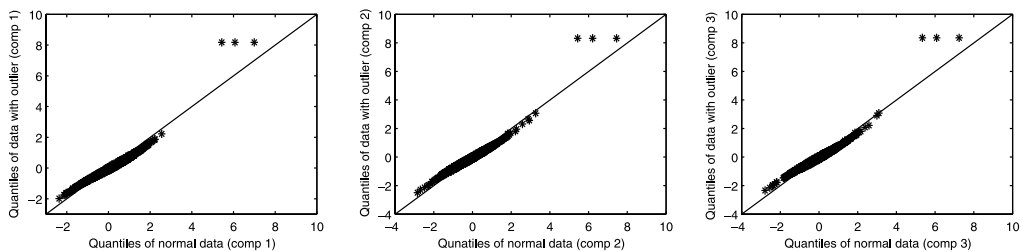}

\caption{The \mbox{\QQ} plots for the two-sample problem, where the first
sample is generated from the standard trivariate
normal distribution, and the second sample contains some outliers.}\label{fig4}
\end{figure}

We next consider an example to demonstrate how our \mbox{\QQ} plots can be
used to detect outliers present in the data. We again consider a
two-sample problem, where the first sample consists of 100 i.i.d.
observations from the standard trivariate normal distribution.
The second sample consists of 97 i.i.d. observations from the standard
trivariate normal distribution and the remaining three data points in
the sample are $(10, 10, 10)$, $(9, 9, 9)$ and $(8, 8, 8)$.
The \mbox{\QQ} plots for this data set are displayed in Figure~\ref{fig4}. The presence
of three outliers in the second sample
is clearly indicated by the plots in Figure~\ref{fig4}.

\subsection{Analysis of real data}\label{sec4.2}

\begin{figure}

\includegraphics{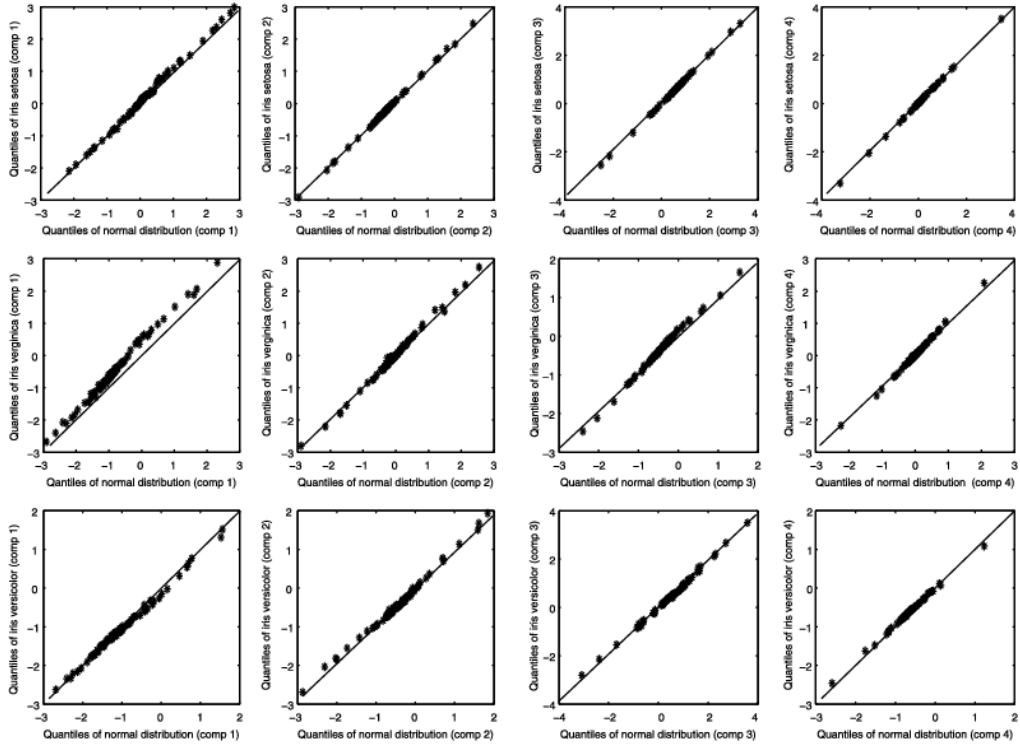}

\caption{The \mbox{\QQ} plots for \emph{Iris setosa} (first row), \emph{Iris
virginica} (second row) and \emph{Iris versicolor} (third row).}\label{fig5}
\end{figure}

We first consider Fisher's \textit{Iris data}, which is available in
\url{http://archive.ics.uci.edu/ml}. In this data, there are three
multivariate samples
corresponding to three different varieties of Iris, namely, \emph{Iris setosa},
\emph{Iris virginica} and \emph{Iris versicolor}. Each sample has size
50. In each sample,
there are four measurements, namely, the sepal length, the sepal width,
the petal length and the
petal width. We would like to determine how close is the distribution
of each sample to a four-dimensional normal distribution. This can be
formulated as a one-sample problem, where $F$ is the distribution of a
sample, and the four-dimensional normal distribution is our specified
distribution $F_0$. Note that $F_0$ involves an unknown mean
$\bolds\mu$ and an unknown dispersion $\Sigma$. For each species,
following the remarks after Theorem \ref{theor2.1}, we estimated
$\bolds\mu$ and $\Sigma$ by the sample mean vector and the sample dispersion
matrix, which are maximum likelihood estimates. Then we standardized
the data in each sample using the corresponding sample mean vector
and the corresponding sample dispersion matrix. The \mbox{\QQ} plots in Figure~\ref{fig5} were constructed using the spatial quantiles of a standardized sample
and the spatial quantiles of the standard four-dimensional normal distribution.

It is visible in the plots in Figure~\ref{fig5} that in almost all cases, the
points are tightly clustered around
the straight line with slope $= 1$ and intercept $= 0$ except in the
first plot for \emph{Iris
virginica}, where the points deviate to some extent from that straight
line. Our one-sample\vadjust{\goodbreak} test for testing $H_0\dvt  F = F_0$ against $H_1\dvt
F\neq F_0$ led to very high $p$-values, namely, $0.841$, $0.413$
and $0.582$ for \emph{Iris setosa}, \emph{Iris virginica} and \emph{Iris
versicolor}, respectively. These $p$-values imply that $H_0$ is to be
accepted, and multivariate normal distributions seem to fit the data
well for all three Iris species.

%
%

\begin{figure}

\includegraphics{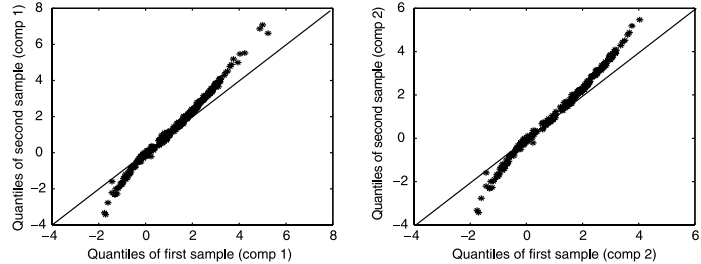}

\caption{The \mbox{\QQ} plots for the vertebral column data.}\label{fig6}
\end{figure}

Our next real data set is the \textit{Vertebral Column data}, which is
available in \url{http://archive.ics.uci.edu/ml/datasets/Vertebral+Column}.
This data set contains six variables on 310 patients, who
belong to two groups. Among the 310 patients, 100 are normal, and the
remaining 210 of them are abnormal.
We view it as a two-sample problem with $F$ as the distribution of the
measurements corresponding to the normal patients, and $G$ as the
distribution of the measurements corresponding to the abnormal patients.
In this study, we considered only two variables, namely, the pelvic
incidence and the pelvic tilt as these two pelvic parameters
are strongly associated with the severity and the stiffness of
lumbosacral spondylolisthesis. Both the pelvic incidence
and the pelvic tilt are angles and measured in the same unit, and no
standardization of the data is necessary in order to compute and
compare the spatial quantiles of these two samples. In Figure~\ref{fig6}, we
display the \mbox{\QQ} plots for
this data. The points in the \mbox{\QQ} plots
are clearly not clustered around any straight line. In fact, most of
the points in each plot lie on a stretched S-shaped curve, which
indicates that the distribution $G$ associated with the abnormal
patients has heavier tails than the distribution $F$ associated with the
normal patients. The $p$-value obtained using the two-sample test for
testing $H_{0}^{*}\dvt  F = G$ against $H_{1}^{*}\dvt  F\neq G$ is $0.038$,
which also indicates that the two distributions are significantly
different.\looseness=-1

%

\begin{figure}

\includegraphics{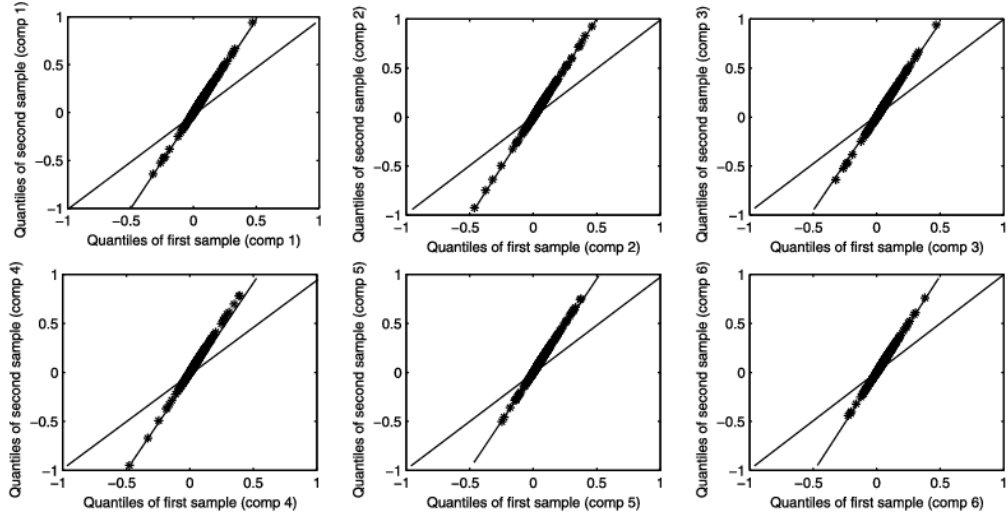}

\caption{The \mbox{\QQ} plots for the monthly sunspot number data.}\label{fig7}
\end{figure}

The third real data set that we consider is the \textit{Monthly Sunspot
number data}, which is available in
\url{http://www.ngdc.noaa.gov/stp/solar/ssndata.html}. This data set contains
monthly average number of sunspots during the period of 1749 to 2009.
As data for 1749 and 2009 are incomplete, we have carried out our analysis
on the observations for the remaining 259 (1750 to 2008) years.
We divided the data into two samples. One sample contains
six-dimensional data corresponding to the six months January, February,
March, October, November and December, and the other one consists
of six-dimensional data corresponding to the months April, May, June,
July, August and September.
The motivation behind splitting
the data into two parts corresponding to the periods October--March and
April--September comes from the fact that one equinox in a year occurs
on March 20--21 and another on September 22--23. We treat this as a
two-sample problem, where $F$ and $G$ are the distributions
corresponding to the sunspot numbers during the periods October--March
and April--September, respectively. The \mbox{\QQ} plots for the data are
presented in Figure~\ref{fig7}.
In each of the plots, the points lie very close to a straight line with
slope $= 2$ and intercept $= 0$. In view of the remark after Theorem
\ref{theor2.2}, these plots indicate that the distributions $F$ and $G$ are
related by the equation $F({\mathbf x}) = G({\mathbf x}/2)$.
Hence, the two multivariate samples corresponding to the two periods
October--March and April--September have distributions that differ only
in the scales
of the variables. The two distributions have the same location, and one
distribution can be obtained from the other by a scale transformation
using the scale factor $2$.
This fact was further confirmed when we carried out some alternative
statistical analysis of the
data such as the comparison of the marginal quantiles and the direct
comparison of the
means and the variances of the variables.

\subsection{Multivariate \mbox{\QQ} plots for data with large dimensions}\label{sec4.3}
When the dimension of the data is large,
there will
be too many two-dimensional plots, and
it will be inconvenient to display and visually examine all of them. In
that case,
one can plot $(l, Q_{{\cal{X}}, l}({\mathbf u}_k) - Q_{{\cal{Y}}, l}({\mathbf
u}_k))$ for $k = 1,\ldots, (n + m)$ and $l = 1,\ldots, d$ in a single
two-dimensional plot with $d$ vertical lines parallel to one another.
We next demonstrate this procedure on some simulated and real data sets.


First, we consider a two-sample problem, where the data in each sample
consists of 10 i.i.d. observations from a standard Brownian motion with
its mean function $m(t) = 0$ and covariance kernel $k(s, t) = \min(s,
t)$, where
$s$, $t\in[0, 1]$ (note that $F = G$ here).
For our second data set, one sample consists of 10 i.i.d. observations
from a standard Brownian motion with its mean function $m_{1}(t) = 0$
and covariance kernel $k_{1}(s, t) = \min(s, t)$ as before (i.e., we
have the same $F$ as before). However, the second sample in the second
data consists of 10 i.i.d. observations from a Brownian motion with its
mean function $m_{2}(t) = 2$ and covariance kernel $k_{2}(s, t) = 2\min
(s, t)$ (which corresponds to the distribution $G$). In our study, we
considered equally spaced points $t_1,\ldots, t_{20}$ in $[0, 1]$ and
sampled the observations at those time points.

The fourth real data set that we consider is the \textit{Sea Level
Pressures data}, which is available in
\url{http://www.cpc.noaa.gov/data/indices/darwin} and
\url{http://www.cpc.noaa.gov/data/indices/tahiti}. This data set consists
of monthly
sea level pressures from two different islands in the southern Pacific ocean,
namely, Darwin ($13^{\circ}$S, $131^{\circ}$E) and Tahiti ($17^{\circ
}$S, $149^{\circ}$W) during
the period 1850--2008. Thus, we have a two-sample problem with each sample
corresponding to an island and containing 159 twelve-dimensional
observations. Here $F$ and $G$ are the distributions of the
multivariate observations corresponding to the two islands. For this
data, each
data point corresponds to a year, and each coordinate of a data point
corresponds to an observation
in a particular month.

\begin{figure}

\includegraphics{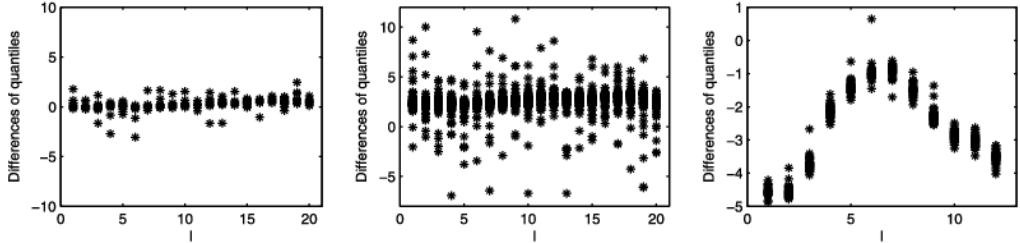}

\caption{The quantile difference plots for the data on Brownian motions
and sea level pressures.}\label{fig8}
\end{figure}

The plots of the quantile differences for the above three data sets are
displayed in Figure~\ref{fig8}. In the first plot
in Figure~\ref{fig8}, the points in each vertical line are tightly clustered
around a horizontal
straight line passing through the origin, which indicates that the
samples are obtained from similar
distributions. It is further confirmed by the large $p$-value $= 0.623$
obtained using our two-sample test
for testing $H_0^{*}\dvt  F = G$ against $H_1^{*}\dvt  F\neq G$. On the other hand,
the difference in the locations and the
scales of the two distributions $F$ and $G$ are clearly visible in the
second plot in Figure~\ref{fig8}. The $p$-value obtained using our two-sample test
in this case is $0.042$, which indicates significant difference between
the two distributions and strong support in favour of $H_1^{*}\dvt  F\neq
G$. It is also amply indicated by the third plot in
Figure~\ref{fig8} as well as the small $p$-value $= 0.045$ obtained using our
two-sample test that the distributions $F$ and $G$ for the two samples
corresponding to the two islands Darwin and Tahiti
are significantly different. 


%


\section{Finite sample level and power study for different tests}\label{sec5}
Here we carry out some simulation studies to compare our tests with the
well-known
multivariate extensions of the Kolmogorov--Smirnov (KS) and the
Cramer--von Mises (CVM) tests (see, e.g., Burke \cite{Bur77} and
Justel, Pe{\~n}a and Zamar \cite{JusPenZam97}) in the one-sample and the two-sample problems.
For testing $H_{0}\dvt  F = F_{0}$ against $H_{1}\dvt  F\neq F_0$, the KS
and the CVM test statistics are $T_{n}^{(1)} =  \sup_{{\mathbf
x}\in\mathbb{R}^{d}}\sqrt{n}|F_{n}({\mathbf x}) - F_{0}({\mathbf x})|$ and
$T_{n}^{(2)} = n \int_{{\mathbf x}\in\mathbb{R}^{d}}[F_{n}({\mathbf
x}) - F_{0}({\mathbf x})]^{2}\mrmd F_{0}({\mathbf x})$, respectively, where
$F_{n}({\mathbf x})$ is the empirical version of $F({\mathbf x})$. To test
$H_{0}^{*}\dvt  F = G$ against $H_{1}^{*}\dvt  F\neq G$, the KS and the CVM
test statistics are $T_{n, m}^{(1)} =  \sup_{{\mathbf x}\in
\mathbb{R}^{d}}\sqrt{n + m}|F_{n}({\mathbf x}) - G_{m}({\mathbf x})|$ and
$T_{n, m}^{(2)} = (n + m) \int_{{\mathbf x}\in\mathbb
{R}^{d}}[F_{n}({\mathbf x}) - G_{m}({\mathbf x})]^{2}\mrmd M_{(n, m)}({\mathbf x})$,
respectively, where
$(n + m)M_{(n, m)}({\mathbf x}) = nF_{n}({\mathbf x}) + mG_{m}({\mathbf x})$, and
$F_{n}$ and $G_{m}$
are the empirical versions of $F$ and $G$, respectively. The KS and the
CVM tests
for multivariate data can be implemented using the asymptotic
distributions of the
corresponding test statistics.

For the one-sample problem, we have considered
$F_{0} = N_{d}$ and $F = (1 - \beta) N_{d} + \beta C_{d}$
and $(1 - \beta) N_{d} + \beta L_{d}$. Here
$\beta\in[0, 1]$, $N_{d}$, $L_{d}$ and $C_{d}$ are the $d$-dimensional
standard normal distribution, the $d$-dimensional Laplace distribution
with p.d.f. $f({\mathbf x}) = (\Gamma(d/2)/2\Gamma(d)\pi^{d/2})\exp
^{-\|{\mathbf x}\|}$ and the $d$-dimensional Cauchy distribution with
p.d.f. $f({\mathbf x}) = (\Gamma((d + 1)/2)/\sqrt{\pi}\Gamma(d/2))(1 +
\|{\mathbf x}\|^{2})^{-(d +1)/2}$, respectively.
In the case of the two-sample problem, we have considered $F = N_{d}$
and $G = (1 - \beta) N_{d} + \beta C_{d}$ and $(1 - \beta) N_{d} +
\beta L_{d}$.

\begin{figure}

\includegraphics{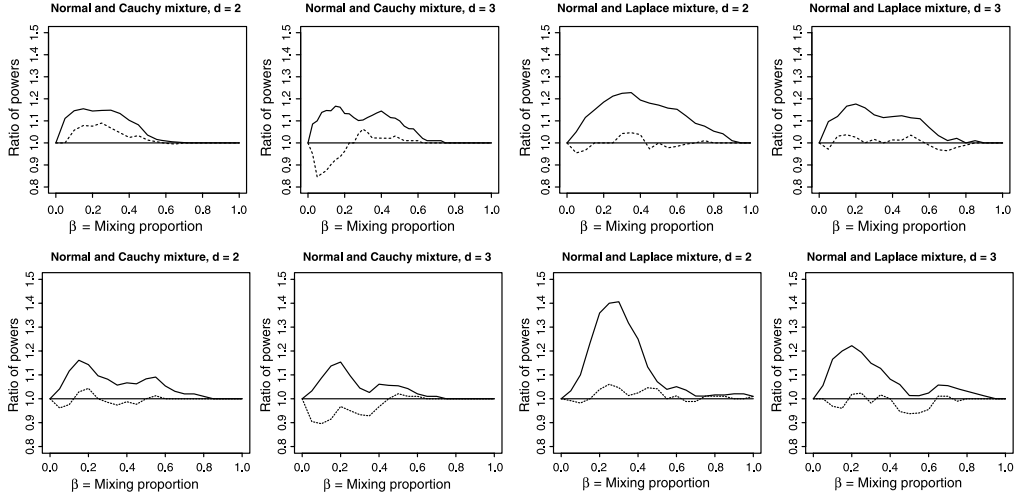}

\caption{The graphs of the ratios of empirical powers based on 1000
Monte Carlo replications at $5\%$ nominal level. The numerator in each
ratio is the power of our test
while the denominators of the ratios corresponding to the \emph{solid}
and the \emph{dotted} curves are the powers
of the KS and the CVM tests, respectively. The first row corresponds to
the one-sample problem with $n = 10$, and the second row corresponds to
the two-sample problem with $n = m = 10$.}\label{fig9}
\end{figure}

In Figure~\ref{fig9}, we have plotted the ratio between the empirical power of
our test (numerator) and that
of another test (denominator) for different values of the parameter
$\beta$.
It is evident from Figure~\ref{fig9} that our test is significantly more
powerful than the KS test in all the cases considered in our simulation
study. However, the CVM test performs better than our test in some
cases, and our test outperforms the CVM test in some other cases.


Friedman and Rafsky \cite{FriRaf79} proposed a multivariate generalization of
the Wald--Wolfowitz run test
using the idea of minimum spanning tree (the MST-run test).
We have compared the empirical powers of our two-sample test with those
of the MST-run test for
$F = N_{d}({\mathbf0}, I_{d})$ and $G = N_{d}(d^{-1/2}\Delta{\mathbf1}_d,
\sigma I_{d})$, where $N_{d}(\bolds\mu,\Sigma)$ is the
$d$-dimensional normal distribution with mean $\bolds\mu$
and dispersion $\Sigma$, ${\mathbf1}_{d}$ is the $d$-dimensional vector of
$1$'s, and the values of $\Delta$ and $\sigma$ are chosen as in
Friedman and Rafsky (\cite{FriRaf81}, page 706).
For sample sizes $n = m = 100$
and $5\%$ nominal level, the results are reported
in Table~\ref{table1}, and it is clear that the MST-run test has inferior performance
compared to our test.

\begin{table}
\caption{Comparison of the empirical powers based on
100 Monte Carlo replications of our two-sample test and the MST-run
test in different dimensions}\label{table1}
\begin{tabular*}{\tablewidth}{@{\extracolsep{\fill}}lllll@{}}
\hline
& $d = 2$ & $d = 5$ & $d = 10$ & $d = 20$\\
& $\Delta= 0.5, \sigma= 1$ & $\Delta= 0.75, \sigma= 1$ & $\Delta=
1.0, \sigma= 1$ & $\Delta= 1.2, \sigma= 1$\\
\hline
Our test & $0.55$ & $0.70$ & $0.83$ & $0.99$\\
MST-run test & $0.35$ & $0.64$ & $0.78$ & $0.86$\\ \hline
& $\Delta= 0, \sigma= 1.2$ & $\Delta= 0, \sigma= 1.2$ & $\Delta=
0, \sigma= 1.1$ & $\Delta= 0, \sigma= 1.075$\\ \hline
Our test & $0.17$ & $0.26$ & $0.07$ & $0.14$\\
MST-run test & $0.14$ & $0.21$ & $0.09$ & $0.13$\\ \hline
\end{tabular*}
\end{table}

For univariate data, our proposed tests in the one-sample and the
two-sample problems lead to new tests that have previously not been
considered in the literature. In addition to the KS and the CVM tests,
there are several other tests that are available in the literature
(see, e.g., Shapiro and Wilk \cite{ShaWil65},
Anderson and Darling \cite{AndDar54} and Ahmad \cite{Ahm93,Ahm96}) for comparing the
distributions of univariate data in the one-sample and the two-sample problems.
We have discussed and compared the performance
of these tests for univariate data in detail in the supplemental
article (see Dhar, Chakraborty and Chaudhuri \cite{DhaChaCha}).

\section{Asymptotic power study under contiguous alternatives}\label{sec6}
Since our tests, the KS and the CVM tests are all
asymptotically consistent, a natural question is how the asymptotic
powers of our tests and the KS
and the CVM tests compare with one another under contiguous
alternatives (see H{\'a}jek and {\v{S}}id{\'a}k \cite{HajSid67}).
%
%
In the case of the one-sample problem, the null hypothesis is
given by
$H_{0}\dvt  F({\mathbf x}) = F_{0}({\mathbf x})$, and we consider a sequence of
contiguous alternatives
$H_{n}\dvt  F({\mathbf x}) = (1 - \gamma/\sqrt{n}) F_{0}({\mathbf x}) + (\gamma
/\sqrt{n}) H({\mathbf x})$
for a fixed $\gamma> 0$ and $n = 1, 2, \ldots\,$.
Consider a multivariate Gaussian process $Z_{1}^{\prime}({\mathbf u})$ with
the mean function
\[
m_{1}({\mathbf u}) = \gamma\bigl[D_{1}\bigl
\{Q_{F_{0}}({\mathbf u})\bigr\}\bigr]^{-1}E_{H} \biggl\{
\frac{{\mathbf x} - Q_{F_{0}}({\mathbf u})}{\|{\mathbf x} - Q_{F_{0}}({\mathbf u})\|} + {\mathbf u} \biggr\}
\]
%
and the covariance kernel $k_{1}({\mathbf u}_1, {\mathbf u}_2)$, where
$k_{1}({\mathbf u}_1, {\mathbf u}_2)$ is as defined before Theorem \ref{theor3.1}. Let
${\cal{V}}^{\prime} = \int\|Z_{1}^{\prime}({\mathbf u})\|^{2}\mrmd {\mathbf u}$, where
the integral is over the same closed ball as in the definition of
$V_{n}$ in Section~\ref{sec3}. We now state a theorem describing the asymptotic
powers of the test based on $V_{n}$ as well as the KS and the CVM tests
under contiguous alternatives.

\begin{theorem}\label{theor6.1}
Assume that $F_{0}$ and $H$ have continuous and
positive densities $f_{0}$ and $h$, respectively, on $\mathbb{R}^{d}
(d\geq2)$, and $E_{F_{0}} \{\frac{h({\mathbf x})}{f_{0}({\mathbf x})} -
1 \}^{2} < \infty$. Then, the sequence of alternatives $H_{n}$
form a contiguous sequence. Under such alternatives, the asymptotic
power of the
test based on $V_{n}$ is given by $P_{\gamma}[{\cal{V}}^{\prime} >
c_{1}(\alpha)]$, where $c_{1}(\alpha)$ is as defined in Theorem
\ref{theor3.1}
such that $P_{\gamma= 0}[{\cal{V}}^{\prime} > c_{1}(\alpha)] = \alpha$.
Further, under those alternatives, the asymptotic powers of the tests
based on $T_{n}^{(1)}$
and $T_{n}^{(2)}$ are given by $P_{\gamma} [ \sup_{{\mathbf
t\in\mathbb{R}^{d}}}|Z_{1}^{\prime\prime}({\mathbf t})| > c_{1}^{*}(\alpha
) ]$
and $P_{\gamma} [ \int_{{\mathbf t}\in\mathbb{R}^{d}}\{
Z_{1}^{\prime\prime}({\mathbf t})\}^{2}\mrmd F_{0}({\mathbf t}) > c_{1}^{**}(\alpha
) ]$, respectively, where $Z_{1}^{\prime\prime}({\mathbf t})$ $({\mathbf
t}\in\mathbb{R}^{d})$ is a Gaussian process with its mean function
$m_{1}^{\prime}({\mathbf t}) = \gamma\{H({\mathbf t}) - F_{0}({\mathbf t})\}$ and
covariance kernel $k_{3}({\mathbf t}_1, {\mathbf t}_2) = F_{0}(\min({\mathbf
t}_1,{\mathbf t}_2)) - F_{0}({\mathbf t}_1)F_{0}({\mathbf t}_2)$. Here
``$\min$'' denotes the coordinatewise minimum of the two vectors in
$\mathbb{R}^{d}$, and $c_{1}^{*}(\alpha)$ and $c_{1}^{**}(\alpha)$
satisfy $P_{\gamma= 0} [ \sup_{{\mathbf t}\in\mathbb
{R}^{d}}|Z_{1}^{\prime\prime}({\mathbf t})| > c_{1}^{*}(\alpha) ] =
\alpha$
and $P_{\gamma= 0} [ \int_{{\mathbf t}\in\mathbb{R}^{d}}\{
Z_{1}^{\prime\prime}({\mathbf t})\}^{2}\mrmd F_{0}({\mathbf t}) > c_{1}^{**}(\alpha
) ] = \alpha$.
\end{theorem}

Next, for the two-sample problem, the null hypothesis is given by
$H_{0}^{*}\dvt  F({\mathbf x}) = G({\mathbf x})$, and we consider a sequence of
alternatives $H_{n, m}^{*}\dvt  G({\mathbf x}) = (1 - \gamma/\sqrt{n + m})
F({\mathbf x}) + (\gamma/\sqrt{n + m}) H({\mathbf x})$ for a fixed $\gamma> 0$
and $n, m = 1, 2, \ldots\,$.
Consider a multivariate Gaussian process
$Z_{2}^{\prime}({\mathbf u})$ with the mean function
\[
m_{2}({\mathbf u}) = -\gamma\bigl[D_{1}\bigl(Q_{F}({
\mathbf u})\bigr)\bigr]^{-1}E_{H} \biggl\{\frac
{{\mathbf y} - Q_{F}({\mathbf u})}{\|{\mathbf y} - Q_{F}({\mathbf u})\|} + {\mathbf
u} \biggr\}
\]
%
and the covariance kernel $k_{2}({\mathbf u}_1, {\mathbf u}_2)$. Here
$k_{2}({\mathbf u}_1, {\mathbf u}_2)$ is as defined before Theorem \ref{theor3.2}. Let
${\cal{T}}^{\prime} = \int\|Z_{2}^{\prime}({\mathbf u})\|^{2}\mrmd {\mathbf u}$, where
the integral is over the same
closed ball as in the definition of $T_{n, m}$ in Section~\ref{sec3}. We now
state a theorem describing the asymptotic powers of the test based on
$T_{n, m}$ as well as the KS and the CVM tests under contiguous alternatives.


\begin{theorem}\label{theor6.2}
Assume that $F$ and $H$ have continuous and
positive densities $f$ and $h$, respectively, on $\mathbb{R}^{d}$
$(d\geq2)$, $E_{F} \{\frac{h({\mathbf y})}{f({\mathbf y})} - 1 \}^{2}
< \infty$, and $n, m\rightarrow\infty$ in such a
way that $ \lim_{n, m\rightarrow\infty}\frac{n}{(n + m)} =
\lambda\in(0, 1)$. Then, the sequence of
densities associated with alternatives $H_{n, m}^{*}$ form a contiguous
sequence. Under such alternatives,
the asymptotic power of the test based on $T_{n, m}$ is given by
$P_{\gamma}[{\cal{T}}^{\prime} > c_{2}(\alpha)]$, where $c_{2}(\alpha)$ is
as defined in Theorem \ref{theor3.2} such that
$P_{\gamma= 0}[{\cal{T}}^{\prime} > c_{2}(\alpha)] = \alpha$.
Further, under those alternatives, the asymptotic powers of the tests
based on $T_{n, m}^{(1)}$ and
$T_{n, m}^{(2)}$ are given by $P_{\gamma} [ \sup_{{\mathbf
t}\in\mathbb{R}^{d}}|Z_{2}^{\prime\prime}({\mathbf t})| > c_{2}^{*}(\alpha
) ]$ and
$P_{\gamma} [ \int_{{\mathbf t}\in\mathbb{R}^{d}}\{
Z_{2}^{\prime\prime}({\mathbf t})\}^{2}\mrmd F({\mathbf t}) > c_{2}^{**}(\alpha
) ]$, respectively, where
$Z_{2}^{\prime\prime}({\mathbf t})$ $({\mathbf t}\in\mathbb{R}^{d})$ is a
Gaussian process
with its mean function $m_{2}^{\prime}({\mathbf t}) = -\gamma\{H({\mathbf t}) -
F({\mathbf t})\}$ and covariance kernel $k_{4}({\mathbf t}_1, {\mathbf t}_2) = \frac
{F(\min({\mathbf t}_1,{\mathbf t}_2)) - F({\mathbf t}_1)F({\mathbf t}_2)}{\lambda(1 -
\lambda)}$.
Here also ``$\min$'' denotes the coordinatewise minimum of the two
vectors in $\mathbb{R}^{d}$, and $c_{2}^{*}(\alpha)$ and
$c_{2}^{**}(\alpha)$ are such that
$P_{\gamma= 0} [ \sup_{{\mathbf t}\in\mathbb
{R}^{d}}|Z_{2}^{\prime\prime}({\mathbf t})| > c_{2}^{*}(\alpha) ] =
\alpha$ and
$P_{\gamma= 0} [ \int_{{\mathbf t}\in\mathbb{R}^{d}}\{
Z_{2}^{\prime\prime}({\mathbf t})\}^{2}\mrmd F({\mathbf t}) > c_{2}^{**}(\alpha
) ] = \alpha$.
\end{theorem}

Theorems \ref{theor6.1} and \ref{theor6.2} enable us to derive the
Pitman efficacies of our tests relative to the KS and the CVM tests.
The Pitman efficacy (see, e.g., Serfling \cite{Ser80} and Lehmann and
Romano~\cite{LehRom05}) of our test relative to another test for
varying choices of the asymptotic power (determined by $\gamma$) is
given by $(\gamma^{\prime}/\gamma)^{2}$, where $\gamma$ and $\gamma
^{\prime}$ are such that the asymptotic power of our test under
contiguous alternatives $(1 - \gamma/\sqrt{n}) F_{0}({\mathbf x}) +
(\gamma/\sqrt{n}) H({\mathbf x})$ (or $(1 - \gamma/\sqrt{n + m})
F({\mathbf x}) + (\gamma/\sqrt{n + m}) H({\mathbf x})$) is the same as
the asymptotic power of the other test under contiguous alternatives
$(1 - \gamma^{\prime}/\sqrt{n}) F_{0}({\mathbf x}) +
(\gamma^{\prime}/\sqrt {n}) H({\mathbf x})$ (or $(1 -
\gamma^{\prime}/\sqrt{n + m}) F({\mathbf x}) + (\gamma^{\prime}/\sqrt
{n + m}) H({\mathbf x})$).

In order to compute the critical values and the powers of our
one-sample and two-sample tests, we have used 1000 simulations of each
Gaussian process and
approximated the integral of the squared norm of a multivariate
Gaussian process by the average
of the squared norms of some appropriate multivariate normal random
vectors. In this numerical study,
we could compute the true covariance matrices as the underlying
distributions were known. We
have computed the critical value and the asymptotic power of the CVM
test in a similar way. However, in the case of the KS test, we have
approximated the supremum of a Gaussian process by a maximum over 1000
simulations of the process.

%
\begin{figure}

\includegraphics{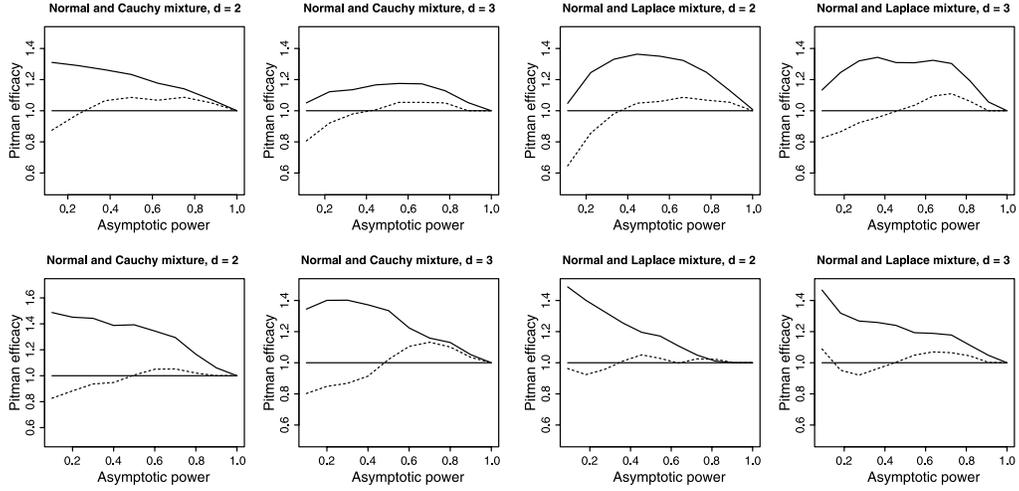}

\caption{The Pitman efficacy of our test relative to the KS test (\emph{solid} curve) and the CVM test (\emph{dotted} curve) at $5\%$ nominal
level. The first row corresponds to the one-sample problem, and the
second row corresponds to the two-sample problem.}\label{fig10}
\end{figure}

In Figure~\ref{fig10}, we have plotted the Pitman efficacy of our test for different
values of the asymptotic power. It is clearly indicated by Figure~\ref{fig10}
that our test
and the CVM test outperform the KS test in terms of the Pitman efficacy
in all the cases
considered here. However, between our test and the CVM test, one has
superior performance in some
cases while the other has superior performance in some other cases, and
there is only a small difference in
their performance.


\begin{appendix}\label{app}
\section*{Appendix: Proofs}\label{sec7}


\begin{pf*}{Proof of Theorem \ref{theor2.1}}
In view of the results in Chaudhuri \cite{Cha96} and Koltchinskii \cite{Kol97},
we have
%
%
\begin{equation}\label{eqA.1}
\sup_{{\mathbf u}} \bigl\|Q_{{\cal{X}}}({\mathbf u}) - Q_{F}({
\mathbf u})\bigr\| = \RMo_{P}(1),
\end{equation}
where the supremum is taken over any given closed ball with the center at
the origin and the radius strictly smaller than one. When $F = F_{0}$,
we have $Q_{F}({\mathbf u}) = Q_{F_{0}}({\mathbf u})$ for all $\|{\mathbf u}\| < 1$.
This along with the uniform convergence result in (\ref{eqA.1}) leads to the
proof of the ``\emph{if} part'' of the theorem.


Next, consider some ${\mathbf u}$ with $\|{\mathbf u}\| < 1$. It follows from
the conditions in the theorem that
with probability tending to one, the spatial rank vectors ${\mathbf
u}_{k}$'s form a dense subset of the unit
ball around the origin as $n\rightarrow\infty$. Since
\[
\lim_{n\rightarrow\infty}P \Biggl( \bigcap_{i = 1}^{d}
\bigl[S_{n,
i}({\cal{X}}, F_{0})\subseteq L(\varepsilon) \bigr]
\Biggr) = 1
\]
for every $\varepsilon> 0$, we must have $Q_{F}({\mathbf u}) = Q_{F_{0}}({\mathbf
u})$ in view of (\ref{eqA.1}).
It now follows from the characterization of multivariate distributions
by the spatial quantiles (see Corollary 2.9 in Koltchinskii (\cite{Kol97}, page
446)) that
$F = F_{0}$. This completes the proof of the ``\emph{only if} part''
of the theorem.
\end{pf*}


\begin{pf*}{Proof of Theorem \ref{theor2.2}}
It follows from the results in Chaudhuri \cite{Cha96} and Koltchinskii \cite{Kol97} that
for the two independent samples ${\cal{X}}$ and ${\cal{Y}}$, we have
$ \sup_{{\mathbf u}}\|(Q_{{\cal{X}}}({\mathbf u}), Q_{{\cal
{Y}}}({\mathbf u})) - (Q_{F}({\mathbf u}), Q_{G}({\mathbf u}))\| = \RMo_{P}(1)$ when
$n$, $m\rightarrow\infty$ in such a way that $ \lim_{n,
m\rightarrow\infty}\frac{n}{(n + m)} = \lambda\in(0, 1)$. Here the
supremum is taken over any given closed ball with the center at the
origin and the
radius strictly smaller than one. Then the proof of the
theorem follows by similar arguments as in the proof of Theorem \ref{theor2.1}.
\end{pf*}

%

\begin{pf*}{Proof of Theorem \ref{theor3.1}}
As proved in Koltchinskii \cite{Kol97}, the centered and normalized
stochastic process $\sqrt{n}\{Q_{{\cal{X}}}({\mathbf u}) -
Q_{F_{0}}({\mathbf u})\}$ converges weakly to the Gaussian process
$Z_{1}({\mathbf u})$ (defined in Section~\ref{sec3}) under $H_{0}$.
Here ${\mathbf u}$ lies in any given closed ball with the center at the
origin and the radius strictly smaller than one. 
It follows from the continuity of the integral functional that $V_{n}$
converges in distribution to ${\cal{V}}$.
Consequently, the asymptotic level of the test will be $\alpha$.

The asymptotic power of the test is given by $ \lim_{n\rightarrow\infty}P_{H_{1}}[V_{n} > c_{1}(\alpha)]$.
Now, note that $V_{n} > c_{1}(\alpha)$ if and only if
$n\int\|\{Q_{{\cal{X}}}({\mathbf u}) - Q_{F_{0}}({\mathbf u})\}-\{Q_{F}({\mathbf
u}) - Q_{F_{0}}({\mathbf u})\}\|^{2}\mrmd {\mathbf u} > c_{1}(\alpha) + n[\int\langle
\{Q_{F}({\mathbf u}) - Q_{F_{0}}({\mathbf u})\}, \{Q_{F}({\mathbf u}) -
Q_{F_{0}}({\mathbf u})\}\rangle \mrmd {\mathbf u} - 2\int\langle\{Q_{{\cal{X}}}({\mathbf
u}) - Q_{F_{0}}({\mathbf u})\}, \{Q_{F}({\mathbf u}) - Q_{F_{0}}({\mathbf u})\}
\rangle \mrmd {\mathbf u}]$.
Here the integrals are over a closed ball with the center at the origin
and the radius strictly smaller than one as before.

When $F\neq F_{0}$, in view of the characterization property of the
spatial quantiles
(see Corollary 2.9 in Koltchinskii \cite{Kol97}), we have $Q_{F}({\mathbf u})\neq
Q_{F_{0}}({\mathbf u})$ for some ${\mathbf u}$
with $\|{\mathbf u}\| < 1$.
The uniform convergence of $Q_{{\cal{X}}}({\mathbf u})$ to $Q_{F}({\mathbf u})$
and the continuity of the spatial quantiles $Q_{F}({\mathbf u})$ and
$Q_{F_{0}}({\mathbf u})$
as functions of ${\mathbf u}$ imply that
$c_{1}(\alpha) + n[\int\langle\{Q_{F}({\mathbf u}) - Q_{F_{0}}({\mathbf u})\},
\{Q_{F}({\mathbf u}) - Q_{F_{0}}({\mathbf u})\}\rangle \mrmd {\mathbf u} - 2\int\langle\{
Q_{{\cal{X}}}({\mathbf u}) - Q_{F_{0}}({\mathbf u})\}, \{Q_{F}({\mathbf u}) -
Q_{F_{0}}({\mathbf u})\}\rangle \mrmd {\mathbf u}]$
tends to $-\infty$ in probability as $n\rightarrow\infty$.
Hence, $P_{H_{1}}[V_{n} > c_{1}(\alpha)]\rightarrow1$ as $n\rightarrow
\infty$. This completes the proof.
\end{pf*}

%


\begin{pf*}{Proof of Theorem \ref{theor3.2}}
Arguing in a similar way as in the proof of Theorem \ref{theor3.1} and
using the weak convergence results in Koltchinskii \cite{Kol97}, and
the independence of the two samples, if $n, m\rightarrow\infty$ in such
a way that $\lambda=  \lim_{n, m\rightarrow\infty}\frac{n}{(n +
m)}\in(0, 1)$, one can show that $T_{n, m}$ converges in distribution
to ${\cal{T}}$ under $H_{0}^{*}$, and consequently, the asymptotic
level of the test that rejects $H_{0}^{*}$ when $T_{n, m} >
c_{2}(\alpha)$ will be $\alpha$. Next, the asymptotic power of the test
is given by $P_{H_{1}^{*}}[T_{n, m} > c_{2}(\alpha)]$. Using similar
arguments as in the second part of the proof of Theorem \ref{theor3.1},
one can establish that $P_{H_{1}^{*}}[T_{n, m} >
c_{2}(\alpha)]\rightarrow1$ as $n, m\rightarrow\infty$.
\end{pf*}

\begin{pf*}{Proof of Theorem \ref{theor6.1}}
The logarithm of the likelihood ratio for testing $H_{0}$ against $H_{n}$
is
%
%
\begin{eqnarray}\label{eqA.2}
L_{n} & = & \sum_{i = 1}^{n}\log
\frac{(1 - \gamma/\sqrt{n})f_{0}({\mathbf
x}_{i}) + (\gamma/\sqrt{n}) h({\mathbf x}_{i})}{f_{0}({\mathbf x}_{i})} = \sum_{i = 1}^{n}\log
\biggl[1 + (\gamma/\sqrt{n}) \biggl\{\frac{h({\mathbf
x}_{i})}{f_{0}({\mathbf x}_{i})} - 1 \biggr\} \biggr]
\nonumber
\\
& = & \frac{\gamma}{\sqrt{n}}\sum_{i = 1}^{n}
\biggl\{\frac{h({\mathbf
x}_{i})}{f_{0}({\mathbf x}_{i})} - 1 \biggr\} - \frac{\gamma^{2}}{2n}\sum
_{i
= 1}^{n} \biggl\{\frac{h({\mathbf x}_{i})}{f_{0}({\mathbf x}_{i})} - 1 \biggr\}
^{2} + R_{n}
\\
& = & \frac{\gamma}{\sqrt{n}}\sum_{i = 1}^{n}k_{i}
- \frac{\gamma
^{2}}{2}\times\frac{1}{n}\sum_{i = 1}^{n}k^{2}_{i}
+ R_{n},\nonumber
\end{eqnarray}
where $k_{i} = \frac{h({\mathbf x}_{i})}{f_{0}({\mathbf x}_{i})} -
1$. Note that $R_{n}\stackrel{P}{\rightarrow} 0$ as $n\rightarrow\infty
$ since $\sigma^{2}:= E_{F_{0}} [\frac{h({\mathbf x})}{f_{0}({\mathbf x})}
- 1 ]^{2} < \infty$. Further, by a straightforward application of
the central limit theorem, the
first term in (\ref{eqA.2}) is asymptotically normal with its mean $= 0$ and
variance $ = \gamma^{2}\sigma^{2}$, and the second term in (\ref{eqA.2})
converges in probability to $\frac{\gamma^{2}}{2}\sigma^{2}$ by the
weak law of large numbers. So, using Slutsky's
theorem, $L_{n}$ is asymptotically normal with mean $= -\frac{\gamma
^{2}}{2}\sigma^{2}$
and variance $= \gamma^{2}\sigma^{2}$. This ensures the contiguity of
the sequence
$H_{n}$ using the corollary to Lecam's first lemma in H{\'a}jek and
{\v{S}}id{\'a}k
(\cite{HajSid67}, pages 204).

Now, we consider ${\mathbf u}_{1},\ldots, {\mathbf u}_{k}$ in a given closed
ball with the center at the origin and the radius strictly
smaller than one, and ${\mathbf t}_{1},\ldots, {\mathbf t}_{l}\in\mathbb
{R}^{d}$. Then, under $H_{0}$, one can establish that
the joint distribution of
$\sqrt{n}\{Q_{{\cal{X}}}({\mathbf u}_{1}) - Q_{F_{0}}({\mathbf u}_{1}),\ldots,
Q_{{\cal{X}}}({\mathbf u}_{k}) - Q_{F_{0}}({\mathbf u}_{k}), F_{n}({\mathbf t}_{1})
- F_{0}({\mathbf t}_{1}),\ldots, F_{n}({\mathbf t}_{l}) - F_{0}({\mathbf t}_{l}),
L_{n}/\sqrt{n}\}$ is asymptotically multivariate normal.
This follows using the Bahadur type linear expansion of $\{Q_{{\cal
{X}}}({\mathbf u}) - Q_{F_{0}}({\mathbf u})\}$
(see Chaudhuri \cite{Cha96}),
the expansion of $L_{n}$ (see (\ref{eqA.2}) above) and the fact that
$F_{n}({\mathbf t}) - F_{0}({\mathbf t})$ is a simple average of i.i.d. random
variables.
Note that for any $p = 1,\ldots, k$, the covariance
between $\sqrt{n}\{Q_{{\cal{X}}}({\mathbf u}_{p}) - Q_{F_{0}}({\mathbf u}_{p})\}
$ and $L_{n}$ is
\begin{eqnarray*}
&&
\frac{\gamma}{n}E_{F_{0}} \Biggl[\sum_{i = 1}^{n}
\biggl\{ D_{1}\bigl[Q_{F_{0}}({\mathbf u}_{p})
\bigr]^{-1} \biggl\{\frac{{\mathbf x}_{i} -
Q_{F_{0}}({\mathbf u}_{p})}{\|{\mathbf x}_{i} - Q_{F_{0}}({\mathbf u}_{p})\|} + {\mathbf u}_{p} \biggr\}
\biggr\}\times \biggl\{\frac{h({\mathbf
x}_{i})}{f_{0}({\mathbf x}_{i})} - 1 \biggr\} \Biggr]
\\
&&\quad
= \gamma\bigl[D_{1}\bigl\{Q_{F_{0}}({\mathbf u}_{p})
\bigr\}\bigr]^{-1}E_{H} \biggl\{\frac
{{\mathbf x} - Q_{F_{0}}({\mathbf u}_{p})}{\|{\mathbf x} - Q_{F_{0}}({\mathbf
u}_{p})\|} + {\mathbf
u}_{p} \biggr\} = m_{1}({\mathbf u}_{p}),
\end{eqnarray*}
because $E_{F_{0}} \{\frac{{\mathbf x} - Q_{F_{0}}({\mathbf u}_{p})}{\|{\mathbf
x} - Q_{F_{0}}({\mathbf u}_{p})\|} + {\mathbf u}_{p} \} = {\mathbf0}$.
Also, one can show that for any $j = 1,\ldots, l$, the covariance
between $\sqrt{n}\{F_{n}({\mathbf t}_{j}) - F_{0}({\mathbf t}_{j})\}$ and
$L_{n}$ is $m_{1}^{\prime}({\mathbf t}_{j}) = \gamma\{H({\mathbf t}_{j}) -
F_{0}({\mathbf t}_{j})\}$.

Now, by a straightforward application of Lecam's third lemma (see
H{\'a}jek and {\v{S}}id{\'a}k \cite{HajSid67}, page 208), one
can establish that under contiguous alternatives, $\sqrt{n}\{Q_{{\cal
{X}}}({\mathbf u}_{1}) - Q_{F_{0}}({\mathbf u}_{1}),\ldots, Q_{{\cal{X}}}({\mathbf
u}_{k}) - Q_{F_{0}}({\mathbf u}_{k})\}$ is asymptotically $kd$-dimensional
multivariate normal with the mean vector having the $d$-dimensional
$p$th block $m_{1}({\mathbf u}_{p})$ ($p = 1, 2,\ldots, k$), and its
$kd\times kd$-dimensional covariance matrix is obtained from the
covariance kernel $k_{1}$, which is given before Theorem \ref{theor3.1}. Further,
the spatial quantile process satisfies the tightness condition under
contiguous alternatives in view of the fact that it is tight under
$H_{0}$. The tightness under $H_{0}$ follows from the weak convergence
of the spatial quantile process (see Koltchinskii \cite{Kol97}).
So, the spatial quantile process $\sqrt{n}\{Q_{{\cal{X}}}({\mathbf u}) -
Q_{F_{0}}({\mathbf u})\}$ converges to $Z_{1}^{\prime}({\mathbf u})$ under
$H_{n}$, where $Z_{1}^{\prime}({\mathbf u})$ is a Gaussian process with its
mean function
$m_{1}({\mathbf u})$ and covariance kernel $k_{1}({\mathbf u}_{1}, {\mathbf
u}_{2})$. Hence, under $H_{n}$, the
asymptotic power of the test based on $V_{n}$ is $P_{\gamma}[{\cal
{V}}^{\prime} > c_{1}(\alpha)]$.

Similarly, using the weak convergence of the stochastic process $\sqrt {n}\{F_{n}({\mathbf t}) - F({\mathbf t})\}$
under $H_{0}$ to a Gaussian process (see, e.g., Bickel and Wichura \cite{BicWic71}) together with Lecam's third lemma,
one can show that under contiguous alternatives, $\sqrt{n}\{F_{n}({\mathbf
t}_{1}) - F_{0}({\mathbf t}_{1}),\ldots, F_{n}({\mathbf t}_{l}) - F_{0}({\mathbf
t}_{l})\}$
is asymptotically $l$-dimensional multivariate normal with the mean
vector having the $j$th component
$m_{1}^{\prime}({\mathbf t}_{j})$ ($j = 1,\ldots, l$), and its $l\times
l$-dimensional covariance matrix is obtained from
the covariance kernel $k_{3}$, which is given in the statement of the theorem.
Now, it follows from the finite-dimensional asymptotic distribution and
the tightness of the process
$\sqrt{n}\{F_{n}({\mathbf t}) - F_{0}({\mathbf t})\}$ under contiguous
alternatives that the stochastic process
$\sqrt{n}\{F_{n}({\mathbf t}) - F_{0}({\mathbf t})\}$ converges
to $Z_{1}^{\prime\prime}({\mathbf t})$ under $H_{n}$, where $Z_{1}^{\prime
\prime}({\mathbf t})$ is a Gaussian process with its mean function
$m_{1}^{\prime}({\mathbf t})$ and covariance kernel $k_{3}({\mathbf t}_1, {\mathbf
t}_2)$. Consequently, under $H_{n}$,
the asymptotic powers of the tests based on $T_{n}^{(1)}$ and
$T_{n}^{(2)}$ are $P_{\gamma} [ \sup_{{\mathbf t}\in\mathbb
{R}^{d}}|Z_{1}^{\prime\prime}({\mathbf t})| > c_{1}^{*}(\alpha) ]$
and $P_{\gamma} [ \int_{{\mathbf t}\in\mathbb{R}^{d}}\{
Z_{1}^{\prime\prime}({\mathbf t})\}^{2}\mrmd F_{0}({\mathbf t}) > c_{1}^{**}(\alpha
) ]$, respectively.
\end{pf*}

\begin{pf*}{Proof of Theorem \ref{theor6.2}}
The logarithm of the likelihood ratio for testing $H_{0}^{*}$ against
$H_{n, m}^{*}$ is
%
%
\begin{eqnarray}\label{eqA.3}
L_{n, m} & = & \log\frac{\prod_{i = 1}^{n}f({\mathbf x}_{i})\prod_{j =
1}^{m}\{(1 - \gamma/\sqrt{n + m})f({\mathbf y}_{j}) + \gamma/\sqrt{n + m}
h({\mathbf y}_{j})\}}{\prod_{i = 1}^{n}f({\mathbf x}_{i})\prod_{j =
1}^{m}f({\mathbf y}_{j})}
\nonumber\\
& = & \sum_{j = 1}^{m}\log \biggl\{1 +
\frac{\gamma}{\sqrt{n + m}} \biggl(\frac{h({\mathbf y}_j)}{f({\mathbf y}_j)} - 1 \biggr) \biggr\}
\\
& = & \frac{\gamma}{\sqrt{n + m}}\sum_{j = 1}^{m}k^{\prime}_{j}
- \frac
{\gamma^{2}}{2(n + m)}\times\sum_{j = 1}^{m}k^{\prime2}_{j}
+ R_{n,
m},\nonumber
\end{eqnarray}
where $k^{\prime}_{j} = \frac{h({\mathbf y}_{j})}{f({\mathbf y}_{j})}
- 1$. Note that $R_{n, m}\stackrel{P}{\rightarrow} 0$ as $n,
m\rightarrow\infty$
since $\sigma_{*}^{2}:= E_{F} \{\frac{h({\mathbf y})}{f({\mathbf y})} -
1 \}^{2} < \infty$. Using similar arguments as in the proof of
Theorem \ref{theor6.1}, $L_{n, m}$ is asymptotically normal with mean $= -\frac
{\gamma^{2}}{2}(1 - \lambda)\sigma_{*}^{2}$ and variance $= \gamma
^{2}(1 - \lambda)\sigma_{*}^{2}$.
This fact ensures the
contiguity of the sequence of densities under $H_{n, m}^{*}$ using the
corollary to Lecam's first lemma in H{\'a}jek and
{\v{S}}id{\'a}k (\cite{HajSid67}, page 204).

Now, here also, we consider ${\mathbf u}_{1},\ldots, {\mathbf u}_{k}$ in a
given closed ball with the center at the origin
and the radius strictly smaller than one, and ${\mathbf t}_{1},\ldots,
{\mathbf t}_{l}\in\mathbb{R}^{d}$. Then, under $H_{0}$,
one can establish that the joint distribution of
$\sqrt{n + m}\{Q_{{\cal{X}}}({\mathbf u}_{1}) - Q_{{\cal{Y}}}({\mathbf u}_{1}),\ldots, Q_{{\cal{X}}}({\mathbf u}_{k}) - Q_{{\cal{Y}}}({\mathbf u}_{k}),
F_{n}({\mathbf t}_{1}) - G_{m}({\mathbf t}_{1}),\ldots, F_{n}({\mathbf t}_{l}) -
G_{m}({\mathbf t}_{l}), L_{n, m}/\sqrt{n + m}\}$
is asymptotically multivariate normal. This asymptotic normality is a
consequence of the independence of the two samples, the
Bahadur type linear expansion of the difference of the spatial
quantiles $Q_{{\cal{X}}}({\mathbf u}) - Q_{{\cal{Y}}}({\mathbf u})$
(see Chaudhuri \cite{Cha96}), 
the expansion of $L_{n, m}$ given in (\ref{eqA.3}) and the
fact that $F_{n}({\mathbf t})$ and $G_{m}({\mathbf t})$ are simple averages of
i.i.d. random variables. Note that for any
$p = 1,\ldots, k$, the covariance between $\sqrt{n + m}\{Q_{{\cal
{X}}}({\mathbf u}_{p}) - Q_{{\cal{Y}}}({\mathbf u}_{p})\}$ and $L_{n, m}$ is
\begin{eqnarray*}
& & E_{F}\Biggl[\sqrt{n + m}\Biggl[\frac{1}{n}\sum
_{i = 1}^{n} \biggl\{D_{1}
\bigl[Q_{F}({\mathbf u}_{p})\bigr]^{-1} \biggl\{
\frac{{\mathbf x}_{i} - Q_{F}({\mathbf u}_{p})}{\|{\mathbf
x}_{i} - Q_{F}({\mathbf u}_{p})\|} + {\mathbf u}_{p} \biggr\} \biggr\}
\\
&&\qquad\hspace*{0pt}{} - \frac{1}{m}\sum_{j = 1}^{m}
\biggl\{D_{1}\bigl[Q_{F}({\mathbf u}_{p})
\bigr]^{-1} \biggl\{\frac{{\mathbf y}_{j} - Q_{F}({\mathbf u}_{p})}{\|{\mathbf
y}_{j} - Q_{F}({\mathbf u}_{p})\|} + {\mathbf u}_{p} \biggr\}
\biggr\}\Biggr] \times\frac{\gamma}{\sqrt{n + m}}\sum_{j = 1}^{m}
\biggl\{\frac{h({\mathbf
y}_{j})}{f({\mathbf y}_{j})} - 1 \biggr\}\Biggr]
\\
&&\quad = -\sqrt{n + m}E_{F}\Biggl[\frac{1}{m}\sum
_{j = 1}^{m} \biggl\{D_{1}
\bigl[Q_{F}({\mathbf u}_{p})\bigr]^{-1} \biggl\{
\frac{{\mathbf y}_{j} - Q_{F}({\mathbf u}_{p})}{\|{\mathbf
y}_{j} - Q_{F}({\mathbf u}_{p})\|} + {\mathbf u}_{p} \biggr\} \biggr\}
\\
&&\hspace*{59.3pt}\qquad{} \times\frac{\gamma}{\sqrt{n + m}}\sum_{j = 1}^{m}
\biggl\{\frac{h({\mathbf
y}_{j})}{f({\mathbf y}_{j})} - 1 \biggr\}\Biggr] \qquad\mbox{(since ${\mathbf x}$ and ${\mathbf
y}$ are independent)}
\\
&&\quad = -\gamma\bigl[D_{1}^{F}\bigl(Q({\mathbf u})\bigr)
\bigr]^{-1}E_{H} \biggl\{\frac{{\mathbf y} -
Q_{F}({\mathbf u}_{p})}{\|{\mathbf y} - Q_{F}({\mathbf u}_{p})\|} + {\mathbf
u}_{p} \biggr\} = m_{2}({\mathbf u}_{p}),
\end{eqnarray*}
because $E_{F} \{\frac{{\mathbf y} - Q_{F}({\mathbf
u}_{p})}{\|{\mathbf y} - Q_{F}({\mathbf u}_{p})\|} + {\mathbf u}_{p} \} = {\mathbf0}$.
Arguing in a similar way as in the proof of Theorem \ref{theor6.1}, one can
establish that under $H_{n, m}^{*}$, the process
$\sqrt{n + m}\{Q_{{\cal{X}}}({\mathbf u}) - Q_{{\cal{Y}}}({\mathbf u})\}$
converges to $Z_{2}^{\prime}({\mathbf u})$,
where $Z_{2}^{\prime}({\mathbf u})$ is a Gaussian process with its mean
function $m_{2}({\mathbf u})$ and covariance
kernel $k_{2}({\mathbf u}_{1}, {\mathbf u}_{2})$, which is defined before
Theorem \ref{theor3.2}. Hence, the asymptotic power of the
test based on $T_{n, m}$ is $P_{\gamma}[{\cal{T}}^{\prime} > c_{2}(\alpha)]$.

Also, under $H_{0}^{*}$, one can show that for any $j = 1,\ldots, l$,
the covariance between
$\sqrt{n + m}\times\{F_{n}({\mathbf t}_{j}) - G_{m}({\mathbf t}_{j})\}$ and $L_{n,
m}$ is
$m_{2}^{\prime}({\mathbf t}_{j}) = -\gamma\{H({\mathbf t}_{j}) - F({\mathbf
t}_{j})\}$. Further,
under $H_{0}^{*}$, the stochastic process $\sqrt{n + m}\{F_{n}({\mathbf t})
- G_{m}({\mathbf t})\}$ converges to a Gaussian process with zero mean and
the covariance kernel
$k_{4}$, which is given in the statement of the theorem (see, e.g.,
Bickel and Wichura~\cite{BicWic71}).
Now, it follows from the finite-dimensional asymptotic distributions
and the tightness of the
process $\sqrt{n + m}\{F_{n}({\mathbf t}) - G_{m}({\mathbf t})\}$ under
contiguous alternatives that the
stochastic process $\sqrt{n + m}\{F_{n}({\mathbf t}) - G_{m}({\mathbf t})\}$
converges to $Z_{2}^{\prime\prime}({\mathbf t})$ under $H_{n, m}^{*}$,
where $Z_{2}^{\prime\prime}({\mathbf t})$ is a Gaussian process with its
mean function $m_{2}^{\prime}({\mathbf t})$
and covariance kernel $k_{4}({\mathbf t}_1, {\mathbf t}_2)$. Consequently,
under $H_{n, m}^{*}$, the asymptotic power of the
test based on $T_{n, m}^{(1)}$ is $P_{\gamma} [ \sup_{{\mathbf t}\in\mathbb{R}^{d}}|Z_{2}^{\prime\prime}({\mathbf t})| >
c_{2}^{*}(\alpha) ]$.

In the case of $T_{n, m}^{(2)}$, we first show that
$(n + m) \int_{{\mathbf x}\in\mathbb{R}^{d}}[F_{n}({\mathbf x}) -
G_{m}({\mathbf x})]^{2}\mrmd (M_{n, m} - F)\stackrel{P}{\rightarrow}0$
as $n, m\rightarrow\infty$ under $H_{0}^{*}$. For that, it is enough to
prove that
$T_{n, m}^{(2, 1)} = (n + m) \int_{{\mathbf x}\in\mathbb
{R}^{d}}[F_{n}({\mathbf x}) - G_{m}({\mathbf x})]^{2}\mrmd (F_{n} - F)\stackrel
{P}{\rightarrow}0$
and $T_{n, m}^{(2, 2)} = (n + m) \int_{{\mathbf x}\in\mathbb
{R}^{d}}[F_{n}({\mathbf x}) - G_{m}({\mathbf x})]^{2}\mrmd (G_{m} - G)\stackrel
{P}{\rightarrow}0$
as $n, m\rightarrow\infty$ under $H_{0}^{*}$. Now, it follows from the
arguments in the proofs of the lemma on page 424 in Kiefer \cite{Kie59} and Theorem 2 in Kiefer and Wolfowitz \cite{KieWol58} that $T_{n, m}^{(2,
1)}\stackrel{P}{\rightarrow}0$ and
$T_{n, m}^{(2, 2)}\stackrel{P}{\rightarrow}0$ as $n, m\rightarrow\infty
$ under $H_{0}^{*}$, and hence,
$(n + m) \int_{{\mathbf x}\in\mathbb{R}^{d}}[F_{n}({\mathbf x}) -
G_{m}({\mathbf x})]^{2}\mrmd (M_{n, m} - F)\stackrel{P}{\rightarrow}0$ as $n,
m\rightarrow\infty$ under $H_{0}^{*}$. Therefore,
$(n + m) \int_{{\mathbf x}\in\mathbb{R}^{d}}[F_{n}({\mathbf x}) -
G_{m}({\mathbf x})]^{2}\mrmd (M_{n, m} - F)\stackrel{P}{\rightarrow}0$ as $n,
m\rightarrow\infty$ under contiguous alternatives $H_{n, m}^{*}$.
Hence, the asymptotic power of the
test based on $T_{n, m}^{(2)}$ under $H_{n, m}^{*}$ is $P_{\gamma}
[ \int_{{\mathbf t}\in\mathbb{R}^{d}}\{Z_{2}^{\prime\prime
}({\mathbf t})\}^{2}\mrmd F({\mathbf t}) > c_{2}^{**}(\alpha) ]$.
\end{pf*}
\end{appendix}

\section*{Acknowledgments}

The research of the first author is partially supported by a grant from
the Council of Scientific and Industrial Research (CSIR), Government of
India. The authors are thankful to two anonymous referees and an
anonymous Associate Editor for several useful comments.

\begin{supplement}
\stitle{Supplement to ``Comparison of multivariate distributions using
quantile--quantile plots
and related tests''}
\slink[doi]{10.3150/13-BEJ530SUPP} 
\sdatatype{.pdf}
\sfilename{BEJ530\_supp.pdf}
\sdescription{In the supplement, we provide additional multivariate
Q--Q plots and discuss the performance of various tests for univariate data.}
\end{supplement}


\printhistory

\end{document}